\documentclass[12pt]{elsarticle}
\usepackage[utf8]{inputenc}

\usepackage{fullpage}

\usepackage{subfig}

\usepackage{amsmath}
\usepackage{amssymb}
\usepackage{latexsym}

\usepackage{hyperref}

\usepackage{boxedminipage}

\usepackage{listings}

\usepackage{ifpdf}

\newcommand{\dx}{\ensuremath{\Delta x}}                     % Delta x
\newcommand{\dt}{\ensuremath{\Delta t}}                     % Delta t
\newcommand{\R}{\ensuremath{\mathbb{R}}}                    % Real field
\newcommand{\wave}{\ensuremath{\mathcal{W}}}
\newcommand{\fwave}{\ensuremath{\mathcal{Z}}}
\newcommand{\cell}{\ensuremath{\mathcal{C}}}
\newcommand{\apdq}{\ensuremath{\mathcal{A}^+ \Delta Q}}
\newcommand{\amdq}{\ensuremath{\mathcal{A}^- \Delta Q}}

\graphicspath{{./}}
\usepackage{epstopdf}

\begin{document}

\ifpdf
\DeclareGraphicsExtensions{.pdf, .jpg, .tif}
\else
\DeclareGraphicsExtensions{.eps, .jpg}
\fi

\begin{frontmatter}

\title{A Numerical Method for the Two Layer Shallow Water Equations with Dry States}

\author[me]{Kyle T. Mandli\fnref{fn1}}
\fntext[fn1]{Present Address: Institute for Computational Engineering and Science, University of Texas at Austin, 201 E 24th ST. Stop C0200, Austin, TX 78712-1229, USA}

\ead{kyle@ices.utexas.edu}
\ead[http://users.ices.utexas.edu/~kyle]{http://users.ices.utexas.edu/~kyle}
\date{\today}
\address[me]{Department of Applied Mathematics, University of Washington, Guggenheim Hall \#414, Box 352420, Seattle, WA 98195-2420, USA}

\begin{abstract}
    A numerical method is proposed for solving the two layer shallow water equations with variable bathymetry in one dimension based on high-resolution f-wave-propagation finite volume methods.  The method splits the jump in the fluxes and source terms into waves propagating away from each grid cell interface.  It addresses the required determination of the system's eigenstructure and a scheme for evaluating the flux and source terms.  It also handles dry states in the system where the bottom layer depth becomes zero, utilizing existing methods for the single layer solution and handling single layer dry states that can exist independently.  Sample results are shown illustrating the method and its handling of dry states including an idealized ocean setting.
\end{abstract}

\end{frontmatter}
% \linenumbers
% ============================================================================
\section{Introduction} \label{sec:introduction} % (fold)

The multilayer shallow water equations have come under increasing interest as a model for primarily long wave phenomena where vertical structure either plays an important role in the flow or the internal structure itself is of interest.  The primary barrier to the use of these equations more broadly has been the complexity and computational cost of the required solvers.  Past approaches to this problem have included the use of more diffusive solvers \cite{Salmon:2002vg}, relaxation approaches \cite{abgrall:2007}, and decoupling schemes that obey entropy laws \cite{Bouchut:2008kn}.  Another difficulty is the loss of hyperbolicity commonly associated with the physical mechanism of Kelvin-Helmholtz instabilities.  Work has also been done to mitigate this effect by applying physically motivated momentum transfer to stabilize the system \cite{CastroDiaz:2010kn}. 

To address some of these difficulties, the numerical method presented here attempts to overcome many of the drawbacks of previous methods while remaining accurate in many scenarios of interest.  In particular this method is intended for oceanic applications where a two-layer model could be used to represent a relatively shallow top ``boundary'' layer and a deeper ``abyssal'' layer.  This implies that the bottom layer can be assumed to become dry before the top layer.  Forcing physics with limited vertical extent such as wind or friction drag are of particular interest, such as in storm surge applications.  In these cases issues such as hyperbolicity are not a concern due to the regime scales being considered.  On the other hand, wetting and drying of the internal surface will occur along the shelf break and must be handled carefully.  Although the step to multiple layers only introduces a limited representation of the three dimensional nature of many of these flows, it is computationally less-expensive than a fully three-dimensional simulation, making the required resolution of some oceanic flows attainable even on modest computing hardware.

The presentation of the numerical method begins with an introduction in section~\ref{sec:multi_layer_shallow_water_equations} to the salient features of the multilayer shallow water equations including methods for evaluating the eigenspace.  A discussion of the basic components of the f-wave approach follows with the specific implementation details for the multilayer shallow water equations and dry states in section~\ref{sec:numerical_approach}.  Finally, example solutions from the numerical method are presented for test cases where dry states are involved in section~\ref{sec:results}.  It should be noted that the presentation will be restricted to the one-dimensional multilayer shallow water equations for clarity.  Many of the salient issues are present in one-dimension and the extension to two-dimensions follows directly from the methods presented here with the formulation of a transverse Riemann solver and extra care dealing with dry-states in the transverse directions.  This topic is left to be presented in future work.

% section introduction (end)
% ============================================================================

% ============================================================================
\section{Multilayer Shallow Water Equations} \label{sec:multi_layer_shallow_water_equations} % (fold)

\begin{figure}[ht]
    \centering
    \includegraphics[width=0.45\textwidth]{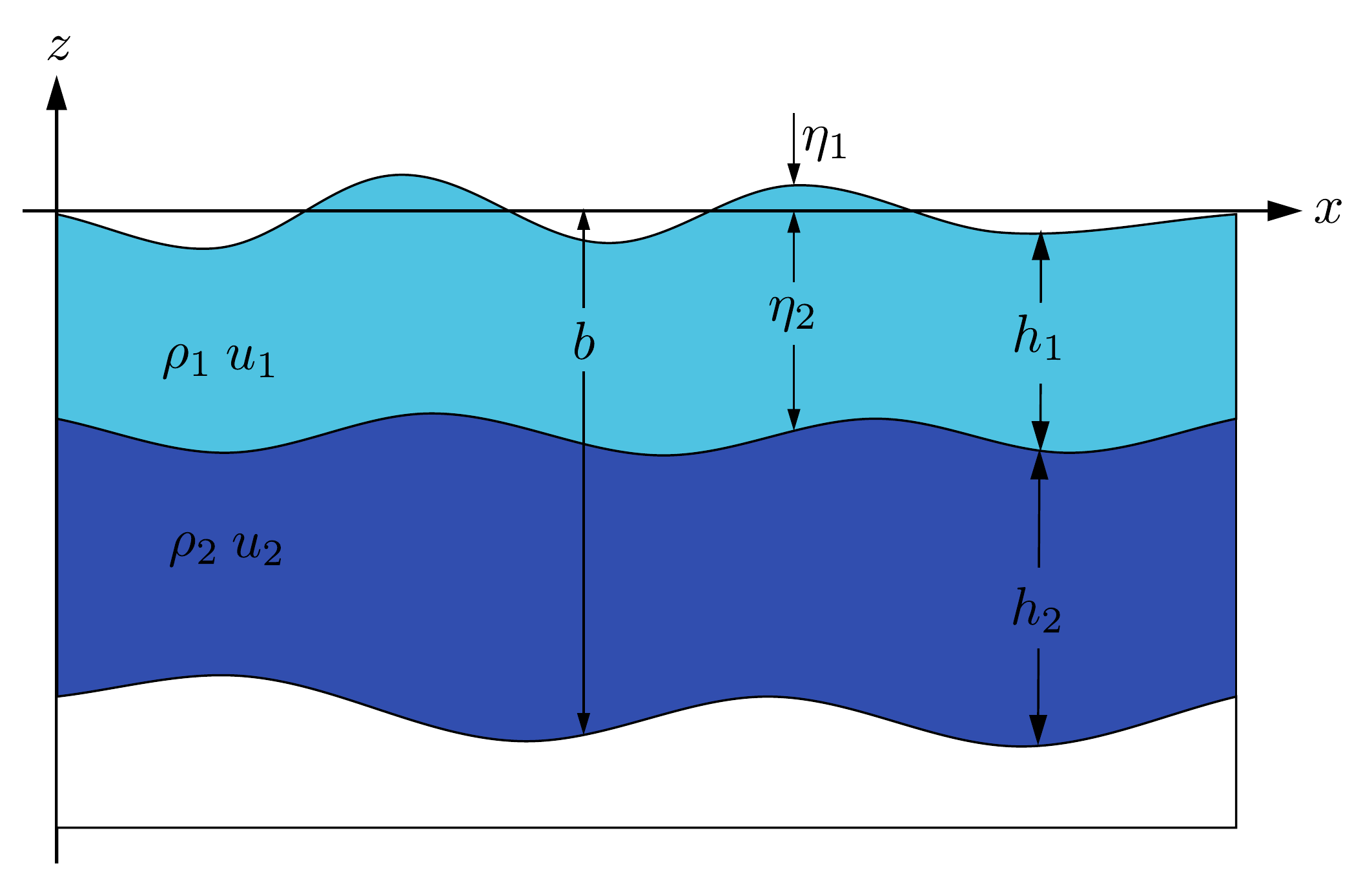}
    \includegraphics[width=0.45\textwidth]{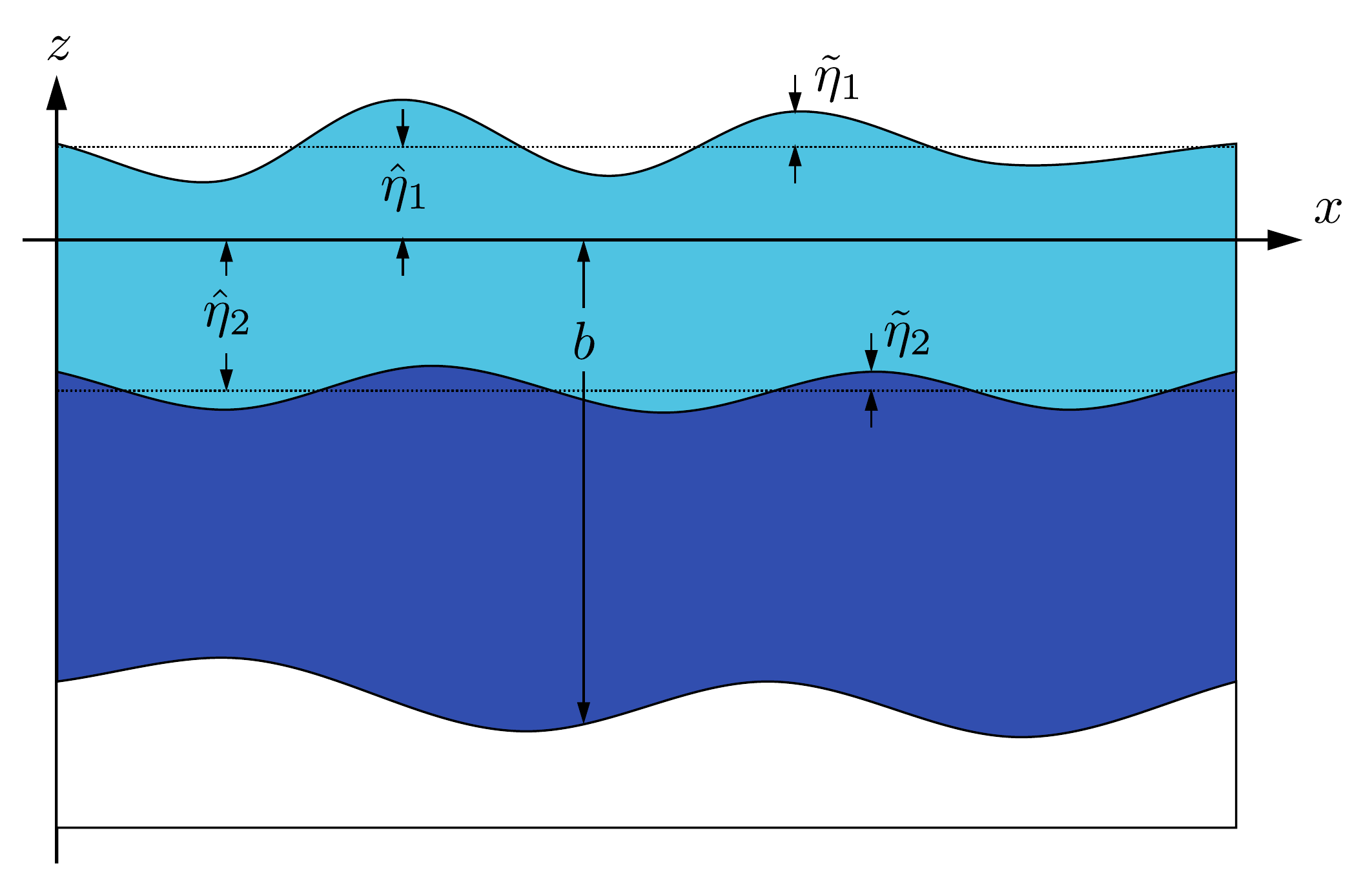} 
    \caption{Coordinates for a one-dimensional system with two-layers and varying bathymetry in general and for the linearized case.}
    \label{fig:multilayer_diagram}
\end{figure}

The multilayer shallow water equations can be derived by integrating the Euler equations in the vertical coordinate direction as in the case of the single-layer equations.  The difference between the multilayer and single-layer shallow water equations is the addition of vertical variation in the density and velocity.  In one-dimension and for two layers the equations are often written as
\begin{equation} \label{eq:mlswe_original}
    \begin{aligned} 
        &(h_1)_t + (h_1 u_1)_x = 0, \\
        &(h_1 u_1)_t + \left(h_1 u_1^2 + \frac{1}{2} g h_1^2 \right)_x = -gh_1(h_2+b)_x, \\
        &(h_2)_t +(h_2 u_2)_x = 0, \\
        &(h_2 u_2)_t + \left(h_2 u_2^2 + \frac{1}{2} g h_2^2 \right)_x = - rgh_2 (h_1)_x - gh_2b_x
    \end{aligned}
\end{equation}
where $h_i$ and $u_i$ are the depths and velocities in each layer respectively, $b$ is the bathymetry from a reference sea-level, $g$ the gravitational acceleration, and $r\equiv \rho_1 / \rho_2$ the ratio of the layer densities (see Figure~\ref{fig:multilayer_diagram}).  Note that we have enumerated the layers with the top layer being indexed first.  The result of the vertical integration and hydrostatic assumption is a system of partial differential equations resembling two sets of single-layer shallow water equations with the addition of a coupling term between the layers.  It is important to note that this coupling is due solely to the integration of the hydrostatic pressure and does not represent momentum transfer due to drag between the layers. 

Another form of equations (\ref{eq:mlswe_original}) involves forgoing the division of the equations by the density of each layer in the derivation and integrating the bottom layer coupling term so that a term appears in the flux of the second-layer's momentum equation rather than as a source term \cite{abgrall:2007}.  In this case, the non-conservative coupling terms in each layer are symmetric.  Using the same notation as before, these equations can be written as
\begin{equation} \label{eq:mlswe_1d}
    \begin{aligned}
        &(\rho_1 h_1)_t + (\rho_1 h_1 u_1)_x = 0,\\
        &(\rho_1 h_1 u_1)_t + \left(\rho_1 h_1 u_1^2 + \frac{1}{2} g \rho_1 h_1^2 \right)_x = -g \rho_1 h_1 (h_2)_x - g \rho_1 h_1 b_x, \\
        &(\rho_2 h_2)_t + (\rho_2 h_2 u_2)_x = 0, \\
        &(\rho_2 h_2 u_2)_t + \left(\rho_2 h_2 u_2^2 + \frac{1}{2} g \rho_2 h_2^2 + g \rho_1 h_2 h_1 \right)_x = g \rho_1 h_1 (h_2)_x - g \rho_2 h_2 b_x.
    \end{aligned}
\end{equation}
The symmetry in the non-conservative products has the benefit that the transfer of momentum due to these coupling terms moves directly between the layers which will be advantageous numerically.  For the remainder of the discussion we will focus our attention solely on the two-layer case solving the system of equations in (\ref{eq:mlswe_1d}).  Extensions of these methods to greater numbers of layers are possible, but for simplicity these complications will be ignored in the rest of the current discussion and is left as future work on the subject.

Since the intention is to consider oceanic applications, an important simplification of the nonlinear equations is the linearization about an ocean at rest.  Taking the steady state where $\hat{u}_1$ and $\hat{u}_2$ are zero and the sea surface $\hat{\eta}_1$ and internal surface $\hat{\eta}_2$ are constant we can rewrite (\ref{eq:mlswe_1d}) as
\begin{equation} \label{eq:mlswe_linearized}
    \begin{aligned}
        &(\tilde{h}_1)_t + (\tilde{\mu}_1)_x = 0, \\
        &(\tilde{\mu}_1)_t + g\hat{h}_1(\tilde{h}_1 + \tilde{h}_2)_x = 0, \\
        &(\tilde{h_2})_t + (\tilde{\mu}_2)_x = 0, ~~~\text{and} \\
        &(\tilde{\mu}_2)_t + g\hat{h}_2 \left[ (\tilde{h}_2)_x + r(\tilde{h}_1)_x \right ] = 0.
    \end{aligned} 
\end{equation}
where we have defined $\hat{h}_1 = \hat{\eta}_1 - \hat{\eta}_2$, $\tilde{h}_1 = \tilde{\eta}_1 - \tilde{\eta}_2$, $\hat{h}_2 = \hat{\eta}_2 - b$, and $\tilde{h}_2 = \tilde{\eta_2}$ for convenience and $\tilde{\mu}_i$ is the perturbation to the momentum of the background ocean at rest such that $\tilde{\mu}_i = \hat{h}_i \tilde{u}_i$ (see figure~\ref{fig:multilayer_diagram}).  Note with these definitions $\hat{h}_2$ is spatially dependent due to the inclusion of $b$.

With these equations, we can write the system \eqref{eq:mlswe_linearized} in the form $\tilde{q}_t + A(\hat{q}) \tilde{q}_x = 0$, where
\[
    \hat{q} = \begin{bmatrix}
        \hat{h}_1 \\ 0 \\ \hat{h}_2 \\ 0
    \end{bmatrix},~~~~\tilde{q} = \begin{bmatrix}
        \tilde{h}_1 \\ \hat{h}_1 \tilde{u}_1 \\ \tilde{h}_2 \\ \hat{h}_2\tilde{u}_2
    \end{bmatrix},~~~\text{and}~~~\tilde{A}(\hat{q})=\begin{bmatrix}
        0 & 1 & 0 & 0 \\
        g\hat{h}_1 & 0 & g \hat{h}_1 & 0 \\
        0 & 0 & 0 & 1 \\
        rg\hat{h}_2 & 0 & g\hat{h}_2 & 0 \\
    \end{bmatrix}.
\]

\subsection{Eigenspace} \label{sub:eigenspace} % (fold)

The eigenspace of hyperbolic PDEs is often of interest and one of the primary sources of difficulties when considering the multilayer shallow water equations.  If one were to directly use the flux Jacobian of (\ref{eq:mlswe_1d}) the eigenvalues and eigenvectors would be identical to two uncoupled shallow water equation systems.  This approach of using a splitting of the layers was shown to be unstable \cite{Castro:2001td} unless suitable corrections are used \cite{Bouchut:2008kn}.  Since the wave speeds predicted by this approach do not take into account the coupling between the layers, this approach is not desirable for methods that depend on this information to construct a Riemann solution.  Instead it is common to write the system in a quasi-linear form $q_t + \tilde{A}(q) q_x = \tilde{S}(q)$ where
\[
    q = [\rho_1 h_1,\rho_1 h_1 u_1,\rho_2 h_2,\rho_2 h_2 u_2]^T
\]
and 
\begin{equation} \label{eq:mlswe_quasilinear_form}
    \tilde{A}(q) = \begin{bmatrix}
        0 & 1 & 0 & 0 \\
        gh_1 - u_1^2 & 2u_1 & rgh_1 & 0  \\
        0 & 0 & 0 & 1\\
        gh_2 & 0 & gh_2 - u_2^2 & 2u_2
    \end{bmatrix} ~~~\text{and}~~~ \tilde{S}(q) = \begin{bmatrix}
        0 \\ -g \rho_1 h_1 b_x \\ 0 \\ - g \rho_2 h_2 b_x
    \end{bmatrix}.
\end{equation}
The characteristic polynomial of the matrix $\tilde{A}(q)$ is then
\begin{equation} \label{eq:mlswe_char_poly}
    ((\lambda-u_1)^2 - gh_1)((\lambda - u_2)^2 - gh_2) - rg^2h_1h_2 = 0.
\end{equation}
Since there exists an expression for the roots of a fourth degree polynomial, it is conceivable that we could evaluate the eigenvalues of the system directly from (\ref{eq:mlswe_char_poly}) as was done in \cite{LAWRENCE:1990ul}.  Numerically this is prohibitively expensive and difficult to do with precision.  Instead, approximations are commonly used to represent the salient features of the system in question and a number of approaches will be presented in section \ref{sub:eigenspace_approximation} including the direct use of LAPACK to calculate the roots. 

Assuming that we can accurately calculate the eigenspeeds $\lambda_p$, we can then find the eigenvectors by solving
% \begin{linenomath}
\begin{equation*}
    \left[ \begin{matrix}
        0 & 1 & 0 & 0 \\
        gh_1 - u_1^2 & 2u_1 & rgh_1 & 0 \\
        0 & 0 & 0 & 1 \\
        gh_2 & 0 & gh_2-u_2^2 & 2u_2
    \end{matrix}\right] \cdot \left [ \begin{matrix}
        1 \\ \alpha_1 \\ \alpha_2 \\ \alpha_3
    \end{matrix}\right] = \left [ \begin{matrix}
        \lambda_p \\ \alpha_1 \lambda_p \\ \alpha_2 \lambda_p \\ \alpha_3 \lambda_p
    \end{matrix}\right].
\end{equation*}
% \end{linenomath}
These equations imply that $\alpha_1 = \lambda_p$ and $\alpha_3 = \lambda_p \alpha_2$.  We then have two equations for one unknown $\alpha_2$ which should satisfy the second and fourth equations simultaneously.  Solving each of the equations separately, we have
\begin{equation} \label{eq:a2_value}
    \alpha_{2,p} = \left\{ \begin{aligned}
        &\frac{(\lambda_p-u_1)^2 - gh_1}{gh_1}, \text{ and} \\
        &\frac{rgh_2}{(\lambda_p-u_2)^2 - gh_2}
    \end{aligned} \right .
\end{equation}
where the subscript of $p$ corresponds to the appropriate eigenvalue.  Setting these two equations equal to one another, we obtain the characteristic polynomial of the system (\ref{eq:mlswe_char_poly}).  We therefore can use either form of $\alpha_{2}$ in equation (\ref{eq:a2_value}) for the eigenvectors.  The final form of the eigenvectors using one of the expressions from (\ref{eq:a2_value}) and dropping the subscript $2$ from $\alpha_{2,p}$ we find
\begin{equation} \label{eq:eigen_vectors_full}
    \left [ 1, \lambda_p, \alpha_p, \lambda_p \alpha_p \right]^T.
\end{equation}

% section multi_layer_shallow_water_equations (end)
% ============================================================================

% ============================================================================
\section{Numerical Approach} \label{sec:numerical_approach} % (fold)

We now describe an algorithm for solving the one dimensional multilayer shallow water equations (\ref{eq:mlswe_1d}).  The method utilizes a f-wave propagation finite volume method, first introduced in \cite{Bale:2002,LeVeque:2001}.  This approach is based on wave-propagation algorithms described in detail in \cite{LeVeque:2002aa} and as such will only be covered briefly here to form a basis for notation and recall the salient ideas.  The rest of this section is dedicated to prescribing the method by which we implement the f-wave propagation method specifically for the multilayer shallow water equations.

\subsection{F-wave Propagation Methods}
\label{sub:f-wave}

Wave-propagation methods are a variant of flux-differencing finite volume methods.  Partitioning the domain into grid cells $\cell$, our goal is to evolve the cell averages $Q^n_i \equiv \frac{1}{\dx_i} \int_{\cell_i} q(x,t^n) dx$ over $\cell_i$ forward in time.  One approach to accomplishing this utilizes the solution of Riemann problems at each grid cell boundary, which are then used to construct numerical fluxes.  The Riemann problems themselves consist of the original system of PDEs on an infinite domain with constant initial data save for a single jump-discontinuity.  The solution to a Riemann problem generally consists of $m$ waves denoted by $\wave^p \in \R^m$ propagating out from the location of the jump-discontinuity at speeds $s^p$.  In the original form of the wave-propagation algorithm these waves are related to the jump-discontinuity at each grid cell interface via
\[
    Q_i - Q_{i-1} = \sum^m_{p=1} \wave^p_{i-1/2}.
\]
The first-order upwind method can then be written as
\[
    Q^{n+1}_{i} = Q^{n}_i - \frac{\dt}{\dx} \left [ \apdq_{i-1/2} - \amdq_{i+1/2} \right]
\]
where $\apdq$ and $\amdq$ represent fluctuations coming from the right and left cells respectively and can be defined in terms of the waves as
\[
    \mathcal{A}^\pm \Delta Q_{i-1/2} = \sum^m_{p=1} (s^p_{i-1/2})^\pm \wave^p_{i-1/2}
\]
where $s^+_{i-1/2} = \text{max}(s^p_{i-1/2},0)$ and $s^-_{i-1/2} = \text{min}(s^p_{i-1/2},0)$.  Extensions to higher-order are possible using limiters applied to each wave such that the update becomes
\[
    Q^{n+1}_{i} = Q^n_i - \frac{\dt}{\dx} (\amdq_{i+1/2} + \apdq_{i-1/2}) - \frac{\dt}{\dx} (\tilde{F}_{i+1/2} - \tilde{F}_{i-1/2})
\]
where 
\[
    \tilde{F}_{i-1/2} = \frac{1}{2} \sum^m_{p=1} |s^p_{i-1/2}| \left (1 - \frac{\dt}{\dx}|s^p_{i-1/2}| \right ) \widetilde{\wave}^p_{i-1/2}
\]
and $\widetilde{\wave}^p_{i-1/2}$ are limited versions of $\wave^p_{i-1/2}$.  If the system being studied is linear, then the $\wave$s are scalar multiples of the eigenvectors with the eigenvalues representing the speeds $s^p$.  In the nonlinear case it is common to use some local linear approximation to the flux Jacobian as is done when Roe-averaging is used.  In either of these cases the waves can be defined as
\[
    \wave^p_{i-1/2} = (\ell^p_{i-1/2})^T(Q_i - Q_{i-1}) r^p_{i-1/2}
\]
where $\ell^p$ and $r^p$ are the left and right eigenvectors of the approximate flux Jacobian $\hat{A}_{i-1/2}$.  This amounts to a projection of the jumps in $Q$ onto the eigenspace of $\hat{A}_{i-1/2}$.  This forms the basis of the wave-propagation method with the primary work being done left to the solution of the Riemann problem at each grid cell boundary. 

An alternative formulation of the wave-propagation method defines the waves in terms of the jump in fluxes instead of the states $q$ such that
\begin{equation} \label{eq:fwave_roe_condition}
    f(Q_i) - f(Q_{i-1}) = \sum^m_{p=1} \fwave^p_{i-1/2}
\end{equation}
where the $\fwave$ are now called f-waves and similar to the $\wave$s we have $\fwave^p \in \R^m$.  In terms of these new waves $\apdq$ and $\amdq$ can be defined as
\[
    \mathcal{A}^\pm \Delta Q_{i-1/2} = \sum^m_{p=1} \text{sgn}(s^p_{i-1/2}) \fwave^p_{i-1/2}.
\]
High-order can again be achieved by applying limiters to each $\fwave$.  As was the case with the $\wave$s, the $\fwave$s can be found by an appropriate linearization of the flux Jacobian matrix.

There are several advantages to using the f-wave formulation over the original wave-propagation formulation.  The first is that the method will be conservative regardless of the linearization of the flux Jacobian employed to calculate the eigenspace.  It also extends naturally to spatially varying flux terms $f(q,x)$.  The advantage that we will make the most use of is the ability to incorporate source terms directly into the waves such that
\[
    f(Q_i) - f(Q_{i-1}) - \tilde{S}_{i-1/2} = \sum^m_{p=1} \fwave^p_{i-1/2}
\]
for some representation of the source term $\tilde{S}$ at $x_{i-1/2}$.

The full f-wave algorithm is implemented identically to the wave-propagation methods except for the differences outlined above.  The basic steps that must be prescribed are:
\begin{enumerate}
    \item State evaluation - The states relevant to the Riemann problem must be extracted from the vectors $Q_\ell$ and $Q_r$.  Also, the determination of the type of dry state problem, if one exists, must be handled.
    \item Compute the eigenvalues $s^p$ and eigenvectors $r^p$.
    \item Compute the jump in the fluxes and source term $F(Q_r) - F(Q_\ell) - \tilde{S} = \delta$.
    \item Project the jump in the fluxes onto the eigenspaces in order to determine the f-waves $\mathcal{Z}^p$.
    \item Calculate the fluctuations $\mathcal{A}^\pm \Delta Q$.
\end{enumerate}

The computation of the eigenspace and flux is dependent on whether a bottom layer dry-state exists and what type it is.  The fully wet case will be presented for each step in its entirety with explanations of appropriate modifications afterwards if a dry-state is present.  The last step of the method, determining the f-waves, fluctuations, and the final update to the grid cell average, is identical to what is presented in \cite{Bale:2002}.

% subsection f_wave_propagation_finite_volume_methods (end)

\subsection{State Evaluation} \label{ssub:state_evaluation} % (fold)
The states that are relevant to the multilayer Riemann problem are the depths $h_i$ of each layer, the momentum $(hu)_i$ in each layer, and the speeds $u_i$ of each layer.  Since the depth and momentum are easily obtained from the state vectors $Q_r$ and $Q_\ell$, these values do not pose a significant challenge.  However, the velocity of each layer needs to be handled carefully since it is obtained by dividing the momentum by the depth.  A simple limiting procedure is therefore used with a dry tolerance $C_{\text{dry}}$ to determine whether to treat the cell as wet or as dry (in which case the velocity is set to zero):
\begin{equation} \label{eq:dry_state_limiter}
    u_{1,2} = \left \{ \begin{aligned}
        &\frac{hu_{1,2}}{h_{1,2}} & &\text{if}~h_{1,2} \ge C_{\text{dry}}& ~~~\text{and} \\
        &0 & &\text{if}~h_{1,2} < C_{\text{dry}}.& \\
    \end{aligned} \right .
\end{equation}
It is important to note that this procedure does not modify the state vectors themselves but only the quantities that will be used to determine the Riemann problem's solution.

In this same stage of the algorithm, the existence and type of dry-states can be determined by also considering the right and left bathymetry states $b_r$ and $b_\ell$.  If a dry state in the bottom layer exists for both sides of the grid cell interface, the problem is treated solely as a single-layer problem and an appropriate single-layer solver is used, including dry states for the single-layer if necessary.  The single-layer solver can handle all instances when the top layer becomes dry due to the assumption that the bottom layer will become dry before the top.

Since the methods described do not provide estimates for the middle states of the Riemann problem only two types of dry-states need to be distinguished: wall dry-states and inundation dry-states.  The difference between these states lies in the relation
\begin{equation} \label{eq:inundation_condition}
    h_{\text{wet}} + b_{\text{wet}} > b_{\text{dry}}
\end{equation}
which represents the case when the fluid on one side of a grid cell interface is high enough that it will flow into the other.  If this is not the case, a wall dry-state exists and no flow of fluid should occur in this layer across the grid cell interface.  Ignoring the middle states is equivalent to ignoring the case where (\ref{eq:inundation_condition}) is not true but the wet cell has enough momentum to overcome the jump in bathymetry.  In this case the cell will have a jump in layer depth greater than the dry cell's bathymetry level and in the subsequent time step, an inundation problem will be solved.
% subsubsection state_evaluation (end)

\subsection{Eigenspace Approximations} \label{sub:eigenspace_approximation} % (fold)

As previously discussed, it is possible that the eigenspace can be evaluated exactly save for the numerical difficulties in solving the characteristic polynomial (\ref{eq:mlswe_char_poly}).  We present here a description of four approaches to approximating the eigenspace that attempts to avoid this difficulty.  These approaches will be compared later in section~\ref{sub:Simple Waves}.

\paragraph{Velocity Difference Expansion} \label{par:velocity_difference_expansion}
One of the most commonly used approximations is based on an expansion about the differences in the layer speeds $u_1 - u_2$ as is done in \cite{Schijf:1953vz}.  Using this approximation and keeping terms in the resulting expansion to first order, we can calculate two sets of eigenspeeds, the external eigenspeeds
\[
    \lambda^{\pm}_{\text{ext}} \approx \frac{h_1u_1 + h_2u_2}{h_1 + h_2} \pm \sqrt{g (h_1 + h_2)}
\]
corresponding to traditional single-layer shallow water waves, also known as barotropic waves, and the internal eigenspeeds
\[
    \lambda^{\pm}_{\text{int}} \approx \frac{h_1u_2 + h_2u_1}{h_1 + h_2} \pm \sqrt{g' \frac{h_1h_2}{h_1+h_2}\left[1-\frac{(u_1 - u_2)^2}{g'(h_1+h_2)}\right]},
\]
where $g' \equiv (1-r)g$, corresponding to waves at the internal surface, often called baroclinic waves.  This approximation has the property that it has a similar form to the eigenvalues of the single-layer equations for the external speeds and the clear result given that the internal wave speeds are small compared to the external speeds when $1-r \ll 1$, \emph{i.e.}, when $\rho_1 \approx \rho_2$.  This approach also leads to the conclusion that if
\[
    \frac{(u_1 - u_2)^2}{ (h_1 + h_2)} \le g'
\]
the approximate eigenvalues will be real, a necessary condition for hyperbolicity.  Assuming that this is true, this also implies that $\lambda^{\pm}_{\text{int}}$ scales as $\sqrt{g'}$ and for the relevant values of $r$ that $\lambda^{\pm}_{\text{int}} \ll \lambda^{\pm}_{\text{ext}}$.

Finding the appropriate states to evaluate the eigenvalues is in general a non-trivial problem.  One possible approach is to allow
\[
    \begin{aligned}
        \lambda^-_{\text{ext,int}} &= \lambda^-_{\text{ext,int}}(Q_\ell)~~\text{and} \\
        \lambda^+_{\text{ext,int}} &= \lambda^+_{\text{ext,int}}(Q_r).
    \end{aligned}
\] 
This approach in general works for systems where the middle states are not expected to differ significantly from the right and left states.  If this were not the case, the estimated speeds are greater than what would be obtained using the true middle states.  This is similar to how Roe averages are found except that in the case of f-waves, the Roe condition 
\[
    s (q_r - q_\ell) = f(q_r) - f(q_\ell)
\]
is automatically satisfied.  Once the eigenvalues are found, the eigenvectors can be evaluated using the same technique, evaluating the eigenvectors corresponding to the left going waves with the left states and similarly with the right states, with the forms found in (\ref{eq:eigen_vectors_full}), which are exact representations once the eigenvalues are known.

\paragraph{Linearized Eigenspace} \label{par:linearized_eigenspace}
Another approach to approximating the eigenvalues is to use the eigenvalues calculated from the characteristic polynomial of the linearized equations (\ref{eq:mlswe_linearized}).  The characteristic polynomial of the linearized system is
\[
    (\lambda^2 - g\hat{h}_1) (\lambda^2 - g\hat{h}_2) - r g^2 \hat{h}_1 \hat{h}_2 = 0.
\]
If we assume that the eigenvectors have the form $v = [1,\lambda,\alpha,\alpha \lambda]^T$ as suggested previously, the relevant equations for the unknowns $\alpha$ and $\lambda$ from $A v - \lambda v = 0$ are
% \begin{linenomath}
\begin{align*}
    &g\hat{h}_1 + g \hat{h}_1 \alpha = \lambda^2 ~~~\text{and} \\
    &rg\hat{h}_2 + g\hat{h}_2 \alpha = \lambda^2 \alpha.
\end{align*}
% \end{linenomath}
We can eliminate $\lambda$ from both equations and find that $\alpha$ must satisfy
\[
    \frac{r\hat{h}_2 + \hat{h}_2 \alpha}{\hat{h}_1 + \hat{h}_1 \alpha} = \alpha.
\]
Dropping the hat notation, we find a quadratic equation for $\alpha$,
\[
    \alpha^2 + \alpha \left(1 - \frac{h_2}{h_1} \right) - r \frac{h_2}{h_1} = 0
\]
leading to
\[
    \alpha_\pm = \frac{1}{2} \left[\gamma - 1 \pm \sqrt{(\gamma-1)^2 + 4r\gamma} \right].
\]
where $\gamma = h_2 / h_1$.  As $\gamma$ approaches 0, implying that the bottom layer depth goes to 0, $\alpha$ takes the value $-1$ or $0$ depending on the wave family in question.  %In another pertinent case, if $\gamma = 1$, implying that the top and bottom layer steady states have equal depth, $\alpha = \pm \sqrt{r}$.

We can now find $\lambda$ in terms of either the top or bottom layers,
\begin{align*}
    \lambda = \left\{ \begin{aligned}
        &\pm \sqrt{gh_1} \sqrt{1+\alpha} ~~~\text{and} \\
        &\pm \sqrt{gh_2} \sqrt{\alpha (r+\alpha)}.
    \end{aligned} \right .
\end{align*}
In order to find which values correspond to which wave-family, consider the dry-state case where the bottom layer goes to zero and take the eigenvalues dependent on the top layer only.  Setting $\gamma=0$ and looking at the positive value of $\alpha$ leads to the eigenvalues $\pm\sqrt{gh_1}$.  This corresponds to the single-layer shallow water wave speeds and can therefore be labeled as the external wave speeds.  If we then examine the case when $\gamma = 1$, then $\alpha = \pm\sqrt{r}$ and $\lambda = \pm \sqrt{gh_1(1 \pm r)}$.  The case where $r > 1$ is eliminated as this would represent a state where the top layer has a higher density than the bottom and would be physically unstable and non-hyperbolic.  In summary, we have
\begin{equation} \label{eq:linearized_evalues}
    \begin{aligned}
        \lambda_{1,\text{ext}} &= -\sqrt{gh_1 (1+\alpha_+)} &~~~~~~&
        \lambda_{2,\text{int}} &= -\sqrt{gh_1 (1+\alpha_-)} \\
        \lambda_{3,\text{int}} &= \sqrt{gh_1 (1+\alpha_-)} &~~~~~~&
        \lambda_{4,\text{ext}} &= \sqrt{gh_1 (1+\alpha_+)}
    \end{aligned}
\end{equation}
where the eigenspeeds have been labeled as described above.  Note that this eigenspace is valid for both the original system we derived (\ref{eq:mlswe_original}) and the modified system (\ref{eq:mlswe_1d}).

The linearized eigensolver uses the same state evaluation as was done for the velocity difference expansion approach.  With the linearization presented in section~\ref{par:linearized_eigenspace} the eigenspace is completely determined by the initial condition (in section~\ref{sec:results} referred to as the linearized static solver), however, and does not change in time since the states used are the steady state values $\hat{h}$ and $\hat{\mu}$.  An alternative to this is to use the full values of the states $h = \hat{h} + \tilde{h}$ and $\mu = \hat{\mu} + \tilde{\mu}$ so that the computed eigenspace evolves in time and is sensitive to changes in the states (in section~\ref{sec:results} referred to as the linearized dynamic solver).  The eigenvectors are evaluated as prescribed by the linearized eigenvectors at the appropriate states as before.

\paragraph{Direct Computation} \label{par:direct_computation}
The final approach we will compare approximates the eigenvalues using a numerical eigensolver such as LAPACK.  The quasi-linear matrix $\tilde{A}(q)$ in (\ref{eq:mlswe_quasilinear_form}) provides the structure for the matrix but the state at which the quasi-linear matrix is evaluated must still be prescribed.  One of the advantages of the f-wave approach is that we do not have to be as careful when choosing this linearization as we are guaranteed that our method is conservative regardless of the linearization used.  In this case the simple arithmetic averages of the state vector $q$ suffices as a means to evaluating the quasi-linear matrix.

\paragraph{Dry States} \label{par:dry_states_eigensolver}
The final issue to resolve in calculating the eigenspace is how to handle dry states.  For the wall dry-state, a modified wall boundary condition is used.  For the bottom layer states, a wall boundary condition that would result in no flux through the step is used, ensuring that the height of the bottom layer used is equal on both sides of the cell interface and the velocity is equal and opposite.  The top layer is affected by the difference between the bathymetry jump and the bottom layer height.  In the inundation case for the bottom layer, the eigenspace is evaluated via the true states in the problem with two possible approaches which differ in what value is used for the dry cell's depth.  If the dry cell's bottom layer depth is set identically to zero, the eigensolver will compute a zero wave speed and the eigenvector will be incorrect.  An estimate to the internal wave is then used based on a single-layer inundation wave,
and 
\[
    s_2 = u_{2r} - 2\sqrt{g(1-r)h_{2r}}
\]
for the case where the left state is dry and
\[
    s_3 = u_{2\ell} + 2 \sqrt{g(1-r)h_{2\ell}}
\]
for the case where the right state is dry.  The eigenvectors are then evaluated using the exact forms (\ref{eq:eigen_vectors_full}) while evaluating with the originally wet state.  The second approach uses non-zero depth for the dry cell's depth to calculate the eigenspace.  In this case it is assumed that the eigenspace is accurately portrayed by a fully wet problem with a small depth, in practice taken to be the dry state cutoff $C_{\text{dry}}$ used earlier.  In both cases, the expectation is that the internal wave will carry the jump in depth in the bottom layer, the external wave moving in the same direction as the dry state is identical to a single-layer shallow water wave for which the eigenspace is known.  

Currently both approaches to inundation appear to produce an entropy violating transonic rarefaction in instances where a strong rarefaction is expected (similar to the all rarefaction solution to the single layer shallow water equations).  A possible course of action would be to modify the appropriate eigenvector that is near $x=0$ by splitting it into two waves but it is currently unclear how to formulate this splitting in the context of an f-wave algorithm.  Fortunately, in the situations of interest these entropy violations do not appear to strongly impact the solution.

% subsection eigenspace_solvers (end)

\subsection{Computation of the Jump in Fluxes} \label{sub:computation_of_the_jump_in_fluxes} % (fold)
If there are no dry states present, the jump in the flux vectors $\delta$ can be computed simply by evaluating
\[
    \delta = f(q_r) - f(q_\ell) - \dx \tilde{S}^n_{i-1/2}
\]
where in the case of the multilayer shallow water equations, the resulting differences are
\begin{equation} \label{eq:delta_jump}
    \delta = \begin{pmatrix}
        [\rho_1 h_1 u_1] \\ 
        [\rho_1 h_1 u_1^2 + 1/2 g \rho_1 h_1^2] + g \rho_1 \overline{h}_1 [h_2 + b] \\
        [\rho_2 h_2 u_2] \\
        [\rho_2 h_2 u_2^2 + 1/2 g \rho_2 h_2^2 + g \rho_1 h_1 h_2] - g \rho_1 \overline{h}_1 [h_2] + g \rho_2 \overline{h}_2 [b]
    \end{pmatrix}
\end{equation}
where $[\cdot]$ and $\overline{\cdot}$ represents the difference and average, respectively, across the cell interface.  These expressions can be derived as path-conservative jump conditions assuming a linear path through state space (see \cite{Mandli:2011te} for this derivation).  

An important property of the differences from \eqref{eq:delta_jump} in conjunction with the f-wave approach is that they lead to a well-balanced algorithm.  In this context, a well-balanced method would be one that would preserve a state at rest where $[\eta_1] = 0$, $[\eta_2] = 0$, $u_1 = 0$, and $u_2 = 0$.  As long as $\delta = 0$, the f-wave algorithm will not perturb the state away from this at rest state.  For $\delta_1$ and $\delta_3$ this condition is trivially satisfied.  For $\delta_2$ and $\delta_4$ it is helpful to consider that for either layer, the general steady state implies that $h_r u_r = h_\ell u_\ell$ and consequently
\[
    [h u^2] = \left [ \frac{(hu)^2}{h}\right] = \frac{|h_r u_r| |h_\ell u_\ell|}{h_r} - \frac{|h_r u_r| |h_\ell u_\ell|}{h_\ell} = -[h] |u_r| |u_\ell|.
\]
It is also useful to observe that
\[
    \frac{1}{2} g [h^2] = \frac{1}{2} g (h_r^2 - h_\ell^2) =  \frac{1}{2} g (h_r + h_\ell) (h_r - h_\ell) = g \overline{h} [h].
\]
Given these expressions, for $\delta_2$ the last three terms are non-trivial but can be rewritten such that
\[
    \delta_2 = -\rho_1 [h_1] |u_{1r}| |u_{1\ell}| + g \rho_1 \overline{h}_1 [h_1] + g \rho_1 \overline{h}_1 ([h_2] + [b])
\]
where the first two terms are zero as $[h_1] = [\eta_1] - [\eta_2] = 0$ and $u_\ell = u_r$ and the last term zero because $[\eta_2] = [h_2] + [b] = 0$.  The expression for $\delta_4$ can rewrite as
\[
    \delta_4 = -\rho_2 [h_2] |u_{2r}| |u_{2\ell}| + g \rho_2 \overline{h}_2 ([h_2] + [b]) + g \rho_1 (\overline{h}_1 [h_2] + \overline{h}_2 [h_1]) - g \rho_1 \overline{h}_1 [h_2]
\]
where using the same relations as before we find that $\delta_4 = 0$.  It should be noted that the single-layer solver employed from \cite{George:2008aa} remains well-balanced for a larger class of steady-states, namely ones where only the assumption that $(hu)_x \equiv 0$ is used, which is not the case for the method presented.

\begin{figure}[tb] %  figure placement: here, top, bottom, or page
    \centering 
    \includegraphics[width=0.3\textwidth]{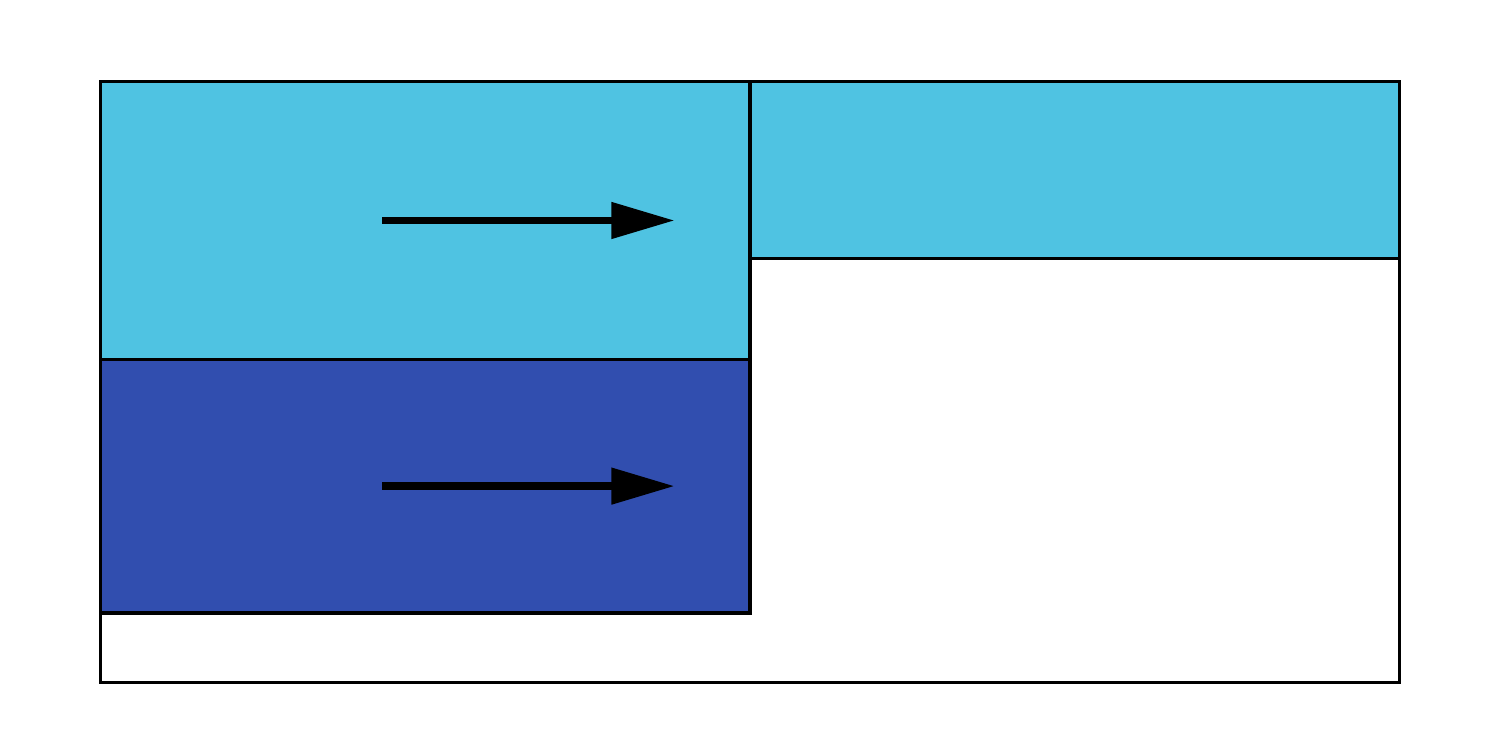}
    \includegraphics[width=0.3\textwidth]{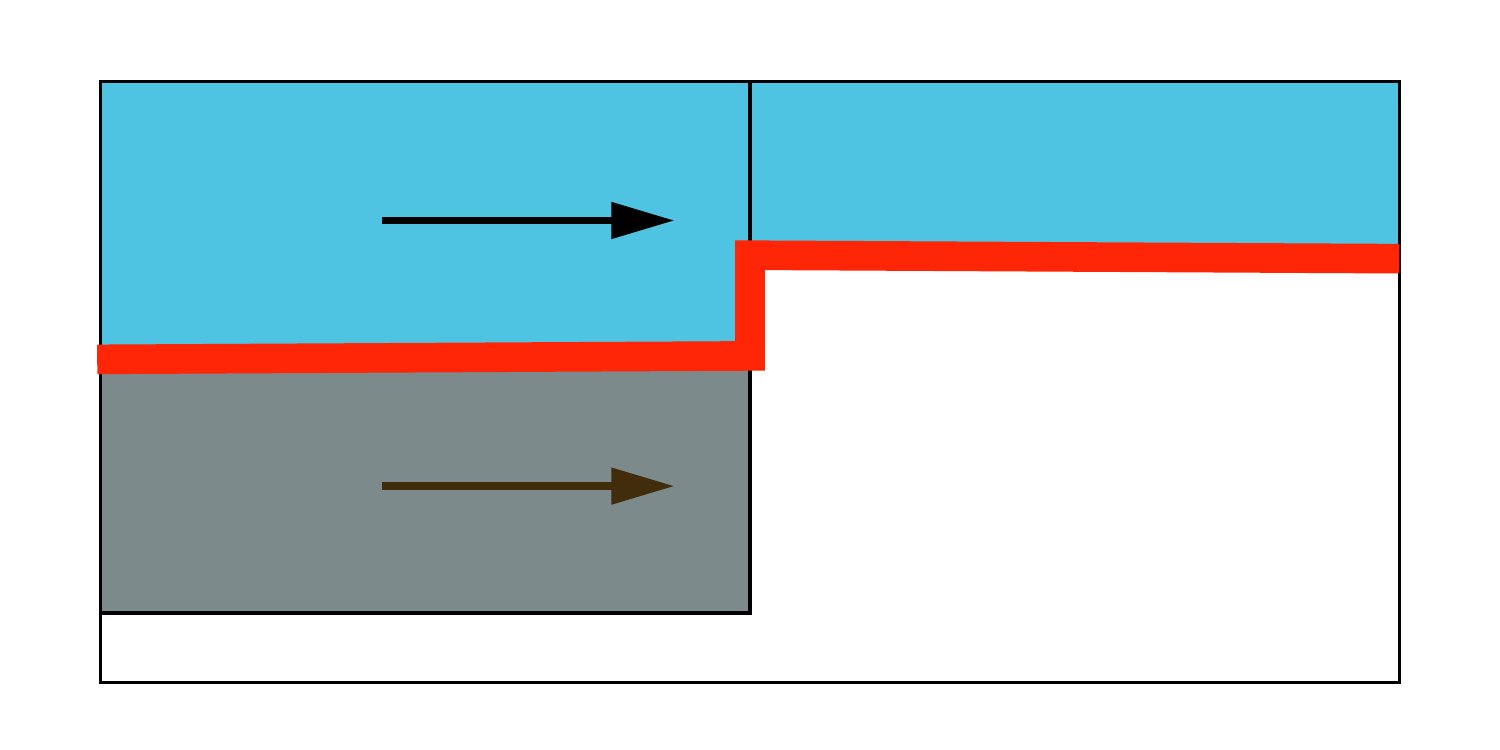}
    \includegraphics[width=0.3\textwidth]{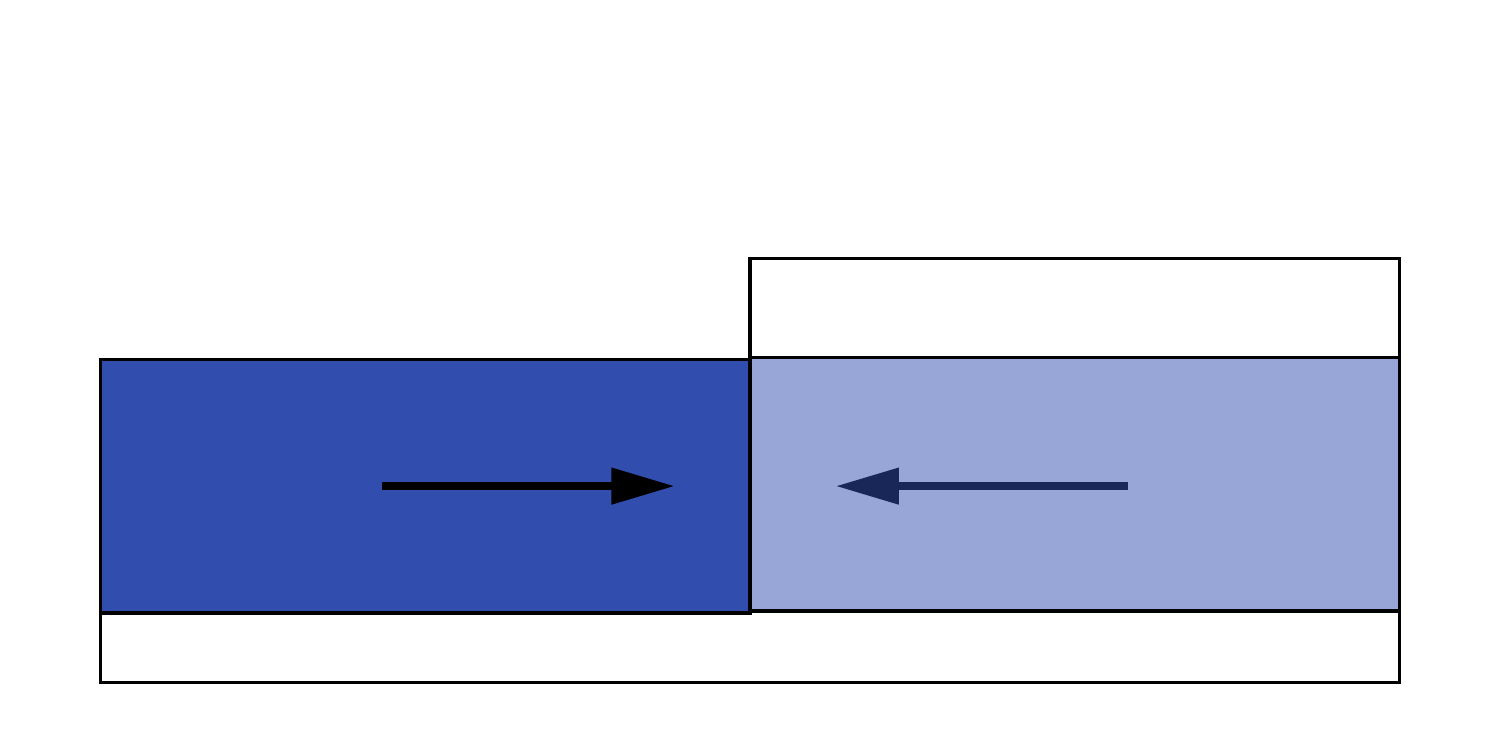}
    \caption{Illustrations of the dry state method presented.  The left most figure represents the general Riemann problem we are interested in solving.  In the middle figure the red line $\eta_2$ is what the top layer ``sees'' when solving the dry state problem and allows for a unified approach to the flux evaluation.  The right most figure is the wall boundary problem solved for the bottom layer.  For the case illustrated a ghost cell with the identical depth and momentum is used to effect a wall.  This also implies that the depth of the top layer $h_1$ remains constant as far as the bottom layer's flux evaluation is concerned.}
    \label{fig:dry_state_method}
\end{figure}

For the wall dry-state problem we consider the case illustrated in figure~\ref{fig:dry_state_method}.  From discussion of the exact Riemann solution to this problem in \cite{Mandli:2011te} we know only some of the states and waves of the solution.  The approximate eigenspace solvers described above try to take this knowledge into account and the primary concern turns to evaluating the jump in the fluxes.

Handling the wall dry state problem is similar to solving the wall dry state problem for the single-layer with some care taken in what states are reflected across the dry state.  Consider the dry-state case illustrated in figure~\ref{fig:dry_state_method}, the goal is to calculate the vector $\delta$ from equation (\ref{eq:delta_jump}) such that the physically relevant Riemann solution is achieved.  In this case, the correct flux jump can be calculated by setting the dry cell's state to one representing a wall boundary condition with
\[
    h_{2r} = h_{2\ell}~~~\text{and} ~~~\mu_{2r} = - \mu_{2r}.
\]
With this state, $\delta_2$ becomes
\[
    \delta_2=\rho_1 [h_1 u_1^2] + 1/2 g \rho_1 [h_1^2] + g \rho_1 \overline{h}_1 (h_{2\ell} - h_{2\ell} + b_r - b_l)
\]
which leads to an incorrect numerical flux.  Looking back at the original modified equations, the top layer is affected by the bottom layer through the momentum exchange term $- g \rho_1 h_1 (h_2)_x$ and combining this with the bathymetry source term $-g\rho_1 h_1 b_x$, the full source term can be rewritten as
\[
    -g \rho_1 h_1 (h_2 + b)_x = -g \rho_1 h_1 (\eta_2)_x.
\]
This implies that the top layer responds to the surface $\eta_2$ only and is insensitive to whether that surface is the bathymetry or the bottom layer.    This allows us to treat the source term by calculating $\eta_2$ with the understanding that if $h_2 = 0$ then $\eta_2 \equiv b$ as illustrated in figure~\ref{fig:dry_state_method}.  Instead of using the wall boundary condition value of $h_{2r}$, the value of $\eta_2$ is calculated since $h_2 + b = \eta_2$.  For this case then the jump can be calculated as
\[
    [h_2 + b] = b_r - h_{2\ell} - b_\ell
\]
and $\delta_2$ is instead calculated by
\[
    \delta_2=\rho_1 [h_1 u_1^2] + 1/2 g \rho_1 [h_1^2] + g \rho_1 \overline{h}_1 (b_r - h_{2\ell} - b_l).
\]

Turning to the momentum flux of the bottom layer, assuming the same boundary values as before $\delta_4$ becomes
% \begin{linenomath}
\begin{align*}
    \delta_4 &= \rho_2 [h_2 u_2^2] + 1/2 g \rho_2 [h_2^2] + \rho_1 [h_1 h_2] - \rho_1 \overline{h}_1 [h_2] + g \rho_2 \overline{h}_2 [b] \\
    &= \rho_1 [h_1 h_2] + g \rho_2 h_{2\ell} (b_r - b_\ell).
\end{align*}
% \end{linenomath}
In this case $\delta_4 \equiv 0$ via the additional condition that $h_{1\ell} = h_{1r}$ and the observation that the bathymetry's effect has been taken into account by the wall boundary condition such that $[b] = 0$.  This is the desired result as there should be no flux through the wall.  The condition $h_{1\ell} = h_{1r}$ is motivated by the observation that the bottom layer is effected only by one side of the Riemann problem and should not be effected by any gradient that may actually be present.

Ensuring that no fluctuations are allowed to wet the bottom layer in the dry grid cell is extremely important.  Small errors below the dry tolerance $C_{\text{dry}}$ can accumulate and may lead to unphysical wetting of dry states.  It may also be necessary to control high-order corrections near large jumps in bathymetry as this can also introduce fluctuations into dry cells that should remain dry.

Another important question is whether the method remains well-balanced in the presence of dry-states in the bottom layer as formulated above.  Since $\delta_4 = 0$ in the wall dry state and inundation would not occur in the steady-state solution considered, the method retains the well-balanced property.
 
% Ensuring in the wall dry state problem that no fluctuations are allowed to wet the bottom layer in the dry grid cell is extremely important.  Small errors that may be below the dry tolerance $C_{\text{dry}}$ will accumulate and for long time simulations may lead to unphysical wetting of dry states.  It also may be necessary to control the high-order correction terms near large jumps in bathymetry as this can introduce small errors in the flux as well.

% subsection computation_of_the_jump_in_fluxes (end)

\subsection{Projection of the Jump in Fluxes onto the Eigenspace} \label{sub:projection_of_the_jump_in_fluxes_onto_the_eigenspace} % (fold)
The projection of the flux difference onto the eigenspace can be obtained by solving
$R \beta = \delta$ where $R$ is the right eigenvector matrix, $\beta$ the resulting projection coefficients, and $\delta$ the flux difference in each field introduced in section~\ref{sub:computation_of_the_jump_in_fluxes}.  The primary difficulty in solving this linear system is ensuring that the solution exists and is unique.  This may not be true numerically if the Riemann problem presented is nearly dry in the bottom layer.  It is important then to identify dry states and treat them appropriately before performing the projection.  In the case of a wall dry-state, the projection collapses down to a simple system of equations since one side of the grid cell is only single-layered and can only carry single-layer shallow water waves.  The inundation case is more difficult and depends on what strategy is used for the eigensolver state evaluation.
\section{Results} \label{sec:results} % (fold)

We now turn to some demonstrative numerical results for the methods outlined above.  The code that has implemented these results can be found at \url{http://www.github.com/clawpack/apps}.  For all the examples variable time steps were taken constrained by the global CFL condition multiplied by 0.9
.  Additionally $C_\text{dry} = 10^{-3}$~meters and $g=9.8~\text{m}/\text{s}^2$ for all examples.  All distances are in meters unless otherwise noted.

\subsection{Simple Waves} % (fold) 
\begin{figure}[htb] %  figure placement: here, top, bottom, or page
    \centering
    \subfloat[Initial condition.]{\includegraphics[width=0.45\textwidth]{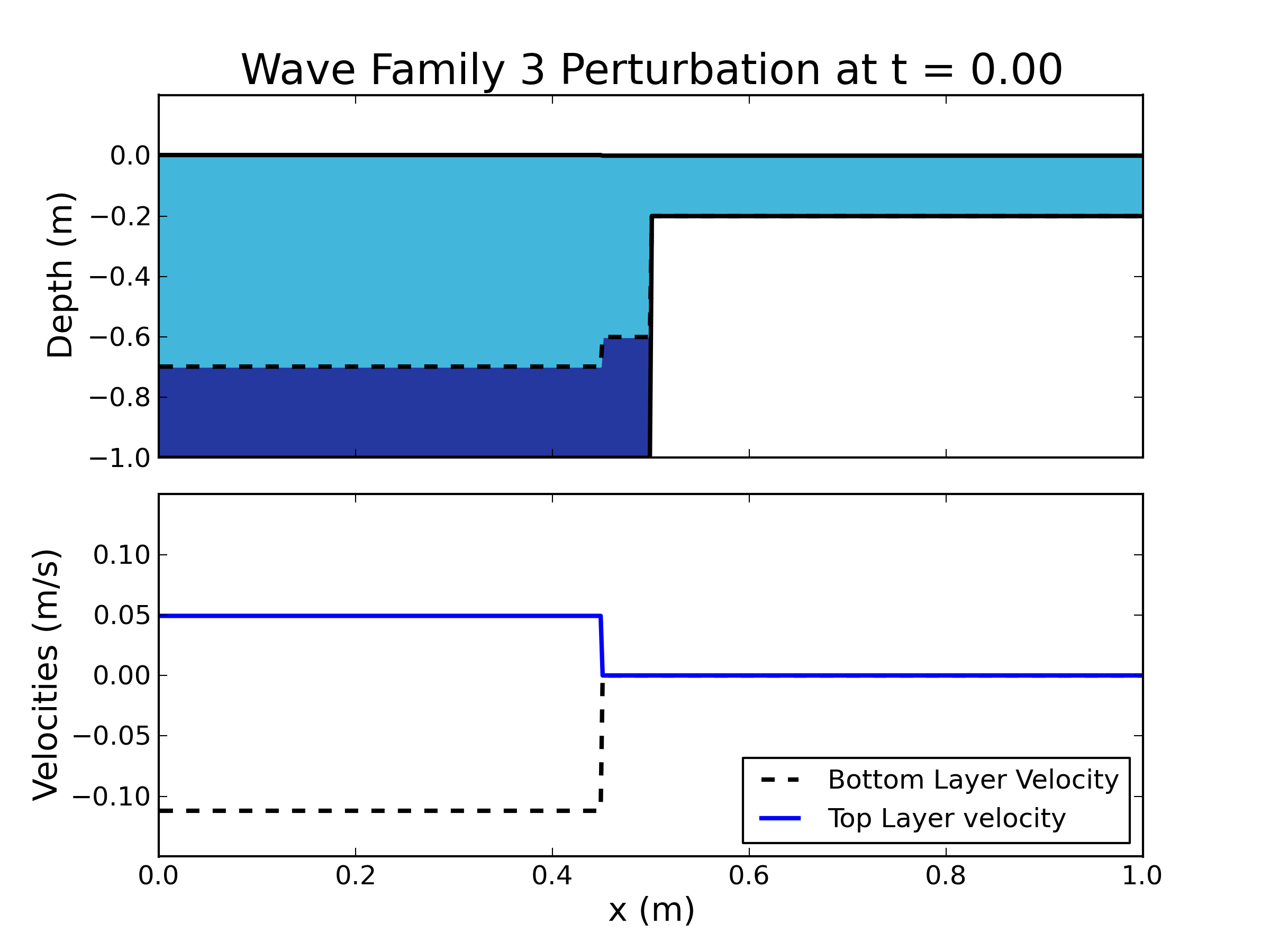} \label{fig:idealized_3_init}}
    \subfloat[]{\includegraphics[width=0.45\textwidth]{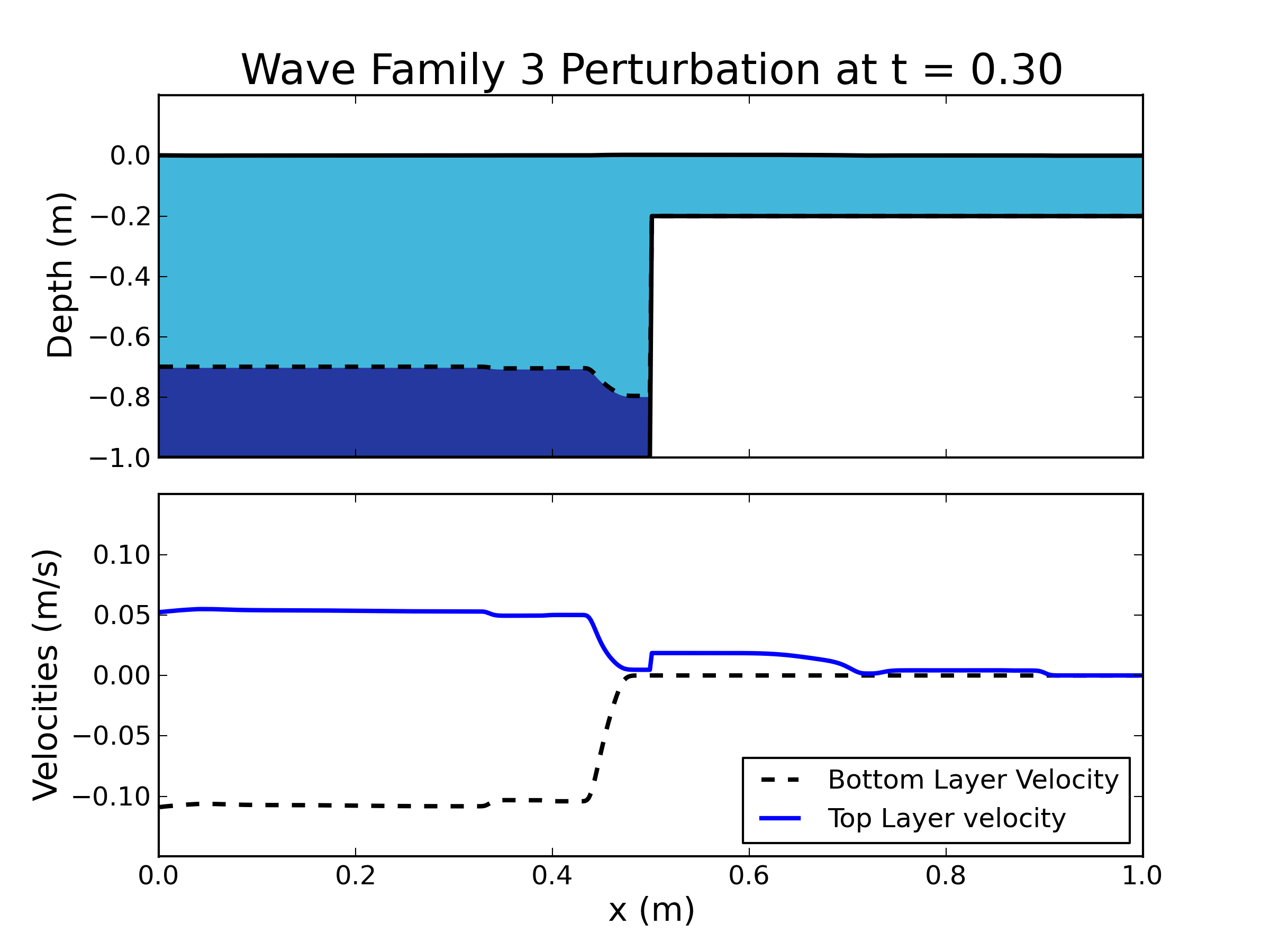} \label{fig:idealized_3_end}}
    \caption{Test showing a perturbation in the 3rd wave-family at a wall dry state problem.  The parameter values used for the test were $r=0.95$ and $\epsilon=0.1$.}
    \label{fig:idealized_3}
\end{figure}

\begin{figure}[htb] %  figure placement: here, top, bottom, or page
    \centering
    \subfloat[Initial condition.]{\includegraphics[width=0.45\textwidth]{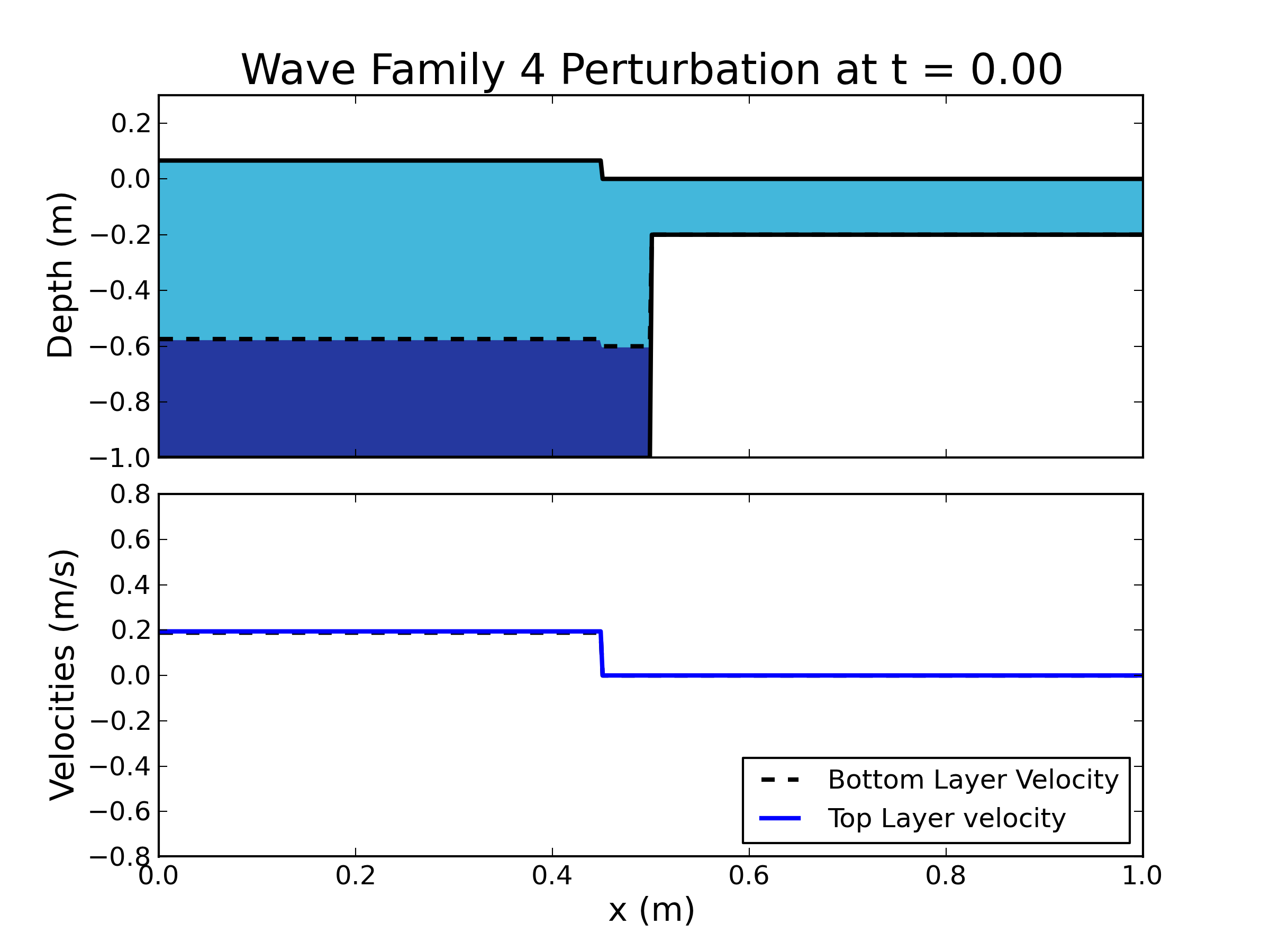}}
    \subfloat[]{\includegraphics[width=0.45\textwidth]{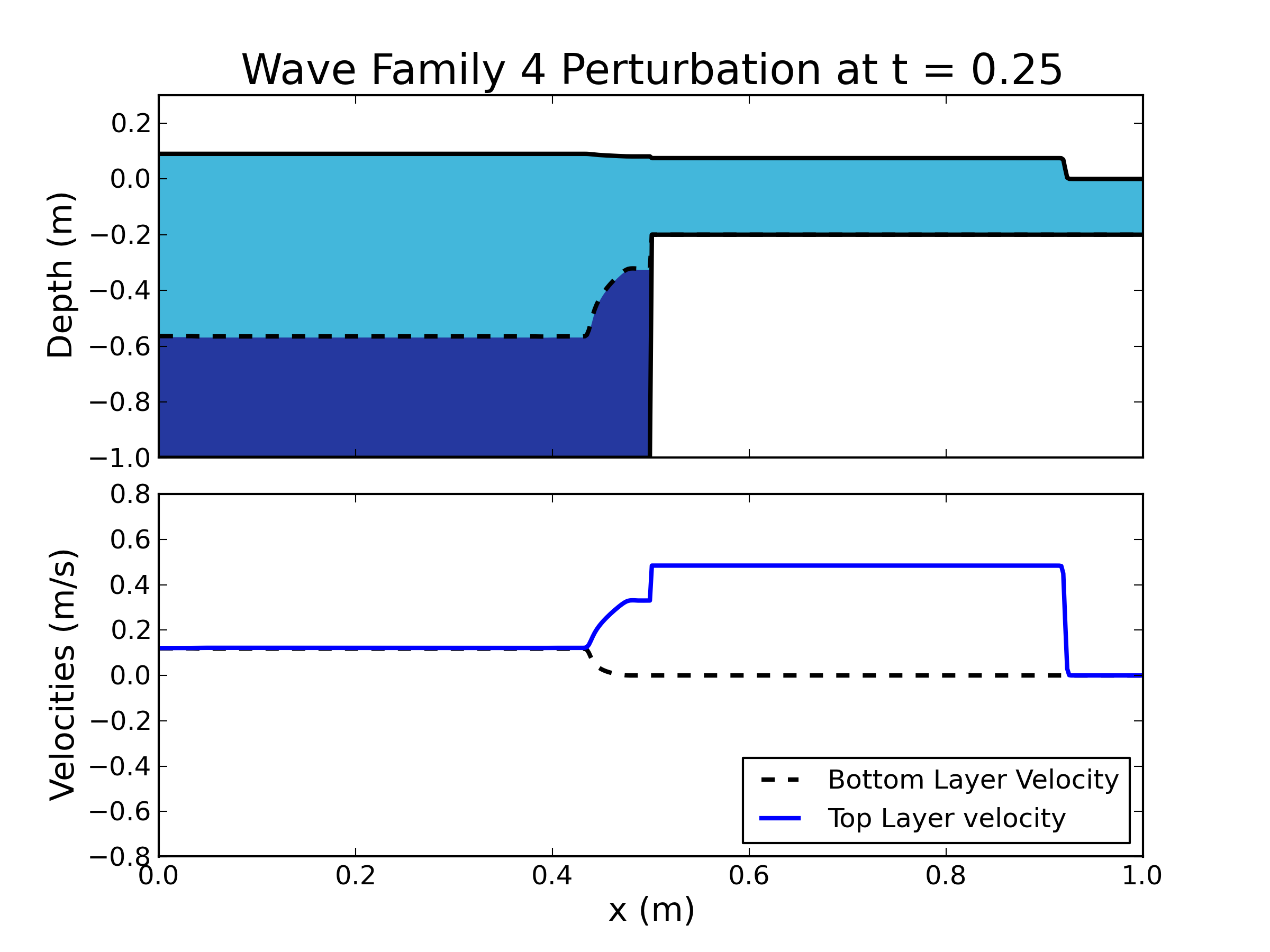}}
    \caption{Test showing a perturbation in the 4th wave-family at a wall dry state problem.  The parameters used for the test were $r=0.95$ and $\epsilon=0.04$.}
    \label{fig:idealized_4}
\end{figure}

\label{sub:Simple Waves}
To illustrate that the method outlined can compute fundamental solutions to the two-layer shallow water equations we calculate the solution in the presence of a dry-state at the center of the domain and an initial condition with a perturbation in only one wave family to a quiescent initial condition $\hat{q}(x,0) = [\rho_1 h_1,0,\rho_2 h_2,0]^T$ such that
\[
    q(x,0) = \left \{ \begin{aligned}
        &r^p\epsilon + \hat{q}(x,0)  &\text{if}~x<x_0& \\
        &\hat{q}(x,0) &\text{if}~x \ge x_0.
    \end{aligned} \right .
\]
These initial conditions are calculated using the approximate linearized eigensystem from section~\ref{par:linearized_eigenspace}.  Figure~\ref{fig:idealized_3} shows a perturbation in the 3rd wave family corresponding to an internal wave with $\hat{\eta}_1 = 0.0$, $\hat{\eta}_2 = -0.6$, $r=0.95$, $x_0=0.45$, $\epsilon=0.1$, $b_r = -0.2$, $b_\ell=-1.0$ and the location of the bathymetry jump at $x=0.5$.  Figure~\ref{fig:idealized_4} shows a perturbation in the 4th wave family representing a wave that should act similar to a single layer shallow water wave with the same parameters as the previous example except for $\epsilon=0.04$.  This smaller value of $\epsilon$ is necessary as a relatively small perturbation can cause a large internal wave that develops after the reflection off of the wall.  For all figures 500 grid cells were used.


\subsubsection{Eigensolver Comparisons and Convergence} \label{ssub:eigen_results} % (fold)

In an effort to show that the methods presented converge, the wave-family tests from section~\ref{sub:Simple Waves} were run with the linearized static, linearized dynamic, velocity difference, and LAPACK eigensolvers at resolutions of 64, 128, 256, 512, and 1024 grid cells and with and without the jump in bathymetry.  Due to the absence of an exact solution to these problems, a ``true'' solution calculated via LAPACK's {\tt DGEEV} routine with 5000 grid cells was used as the basis for comparison.

\begin{figure}[thb]
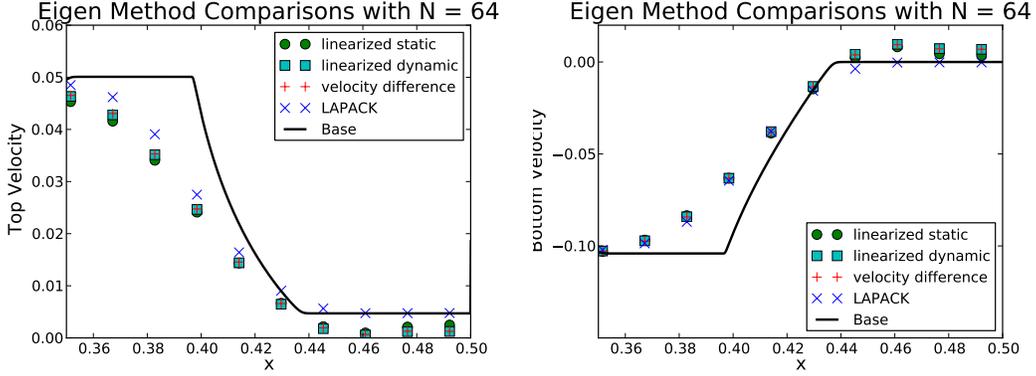

    % \begin{center}
    \centering
        \includegraphics[width=0.42\textwidth]{fig_5_u_top_n64.png}
        \includegraphics[width=0.42\textwidth]{fig_5_u_bottom_n64.png}
    % \end{center}
    \caption{Zoomed-in section of the 3rd wave family dry-state test case with 64 grid cells.  The figures show the velocity computed by each method of the top and bottom layers respectively.}
    \label{fig:eigen_method_comparison}
\end{figure}

Convergence of all the methods behaved nearly identically in all test cases.  Observed convergence rates for the test cases where smooth solutions exist were found to be second order or better when convergence was expected (\emph{i.e.}, when a discontinuity was not present in the solution) and nearly first order otherwise (see tables~\ref{table:convergence_wet_wave} and~\ref{table:convergence_dry_wave}).  The examples demonstrate the ability of the method to converge to traveling simple wave solutions in each of the families with and without a dry-state present.  The order was calculated by fitting a linear polynomial through the given data points minimizing the squared error.  The higher orders of convergence for the 4th wave family should probably not be taken as an indication of the method's convergence order but rather a consequence of the choice of a ``true'' solution.  The 4th wave family convergence should match a single-layer shallow water approach but its order of convergence is higher than expected.  The largest discrepancy between methods can be seen in the velocity comparisons where the linearized and velocity difference approaches yielded slightly different velocity fields than the LAPACK eigensolver as in figure~\ref{fig:eigen_method_comparison}.  The convergence order in the 3rd wave family is expected as the perturbation is given to the internal surface where the bathymetry jump reflects and a first-order error may be present.  In this analysis it should be noted that although all the methods compared revert to the linearized eigensolvers at dry-state interfaces in the bottom layer, the velocity difference eigensolver failed in the 4th wave-family dry state problem after the wave arrives at the wall.  
                                                       
\begin{table}[tb]
    \begin{center}
        \begin{tabular}{||l|l|r|r||}
            \hline \hline 
            ~ & ~ & \multicolumn{2}{c||}{Order} \\
            Field                   &                Method &   3rd Family & 4th Family \\ 
            \hline \hline
            Top Layer Depths        & linearized static     &   1.65 &   3.18 \\ 
                                    & linearized dynamic    &   1.65 &   3.18 \\ 
                                    & velocity difference   &   1.47 &   3.23 \\ 
                                    & LAPACK {\tt DGEEV}    &   1.58 &   3.24 \\ 
            \hline 
            Top Layer Velocities    & linearized static     &   1.67 &   4.06 \\ 
                                    & linearized dynamic    &   1.67 &   4.06 \\ 
                                    & velocity difference   &   1.49 &   3.91 \\ 
                                    & LAPACK {\tt DGEEV}    &   1.60 &   3.70 \\ 
            \hline 
            Bottom Layer Depths     & linearized static     &   2.30 &   1.59 \\ 
                                    & linearized dynamic    &   2.30 &   1.59 \\ 
                                    & velocity difference   &   2.15 &   1.64 \\ 
                                    & LAPACK {\tt DGEEV}    &   2.24 &   1.57 \\ 
            \hline 
            Bottom Layer Velocities & linearized static     &   1.74 &   1.59 \\ 
                                    & linearized dynamic    &   1.74 &   1.59 \\ 
                                    & velocity difference   &   1.56 &   1.64 \\ 
                                    & LAPACK {\tt DGEEV}    &   1.66 &   1.61 \\ 
            \hline \hline
        \end{tabular}
    \end{center}
    \caption{Calculated order of test cases run without a jump in bathymetry.
    \label{table:convergence_wet_wave}}
\end{table}

\begin{table}[thb]
    \begin{center}
        \begin{tabular}{||l|l|r|r||}
            \hline \hline 
            ~ & ~ & \multicolumn{2}{c||}{Order} \\
            Field                   &                Method &  3rd Family & 4th Family  \\ 
            \hline \hline 
            Top Layer Depths        & linearized static     &   1.49 &   2.24 \\ 
                                    & linearized dynamic    &   1.49 &   2.24 \\ 
                                    & velocity difference   &   1.25 &   2.37 \\ 
                                    & LAPACK {\tt DGEEV}    &   1.30 &   2.36 \\ 
            \hline 
            Top Layer Velocities    & linearized static     &   1.57 &   2.25 \\ 
                                    & linearized dynamic    &   1.57 &   2.25 \\ 
                                    & velocity difference   &   1.32 &   2.41 \\ 
                                    & LAPACK {\tt DGEEV}    &   1.37 &   2.40 \\ 
            \hline 
            Bottom Layer Depths     & linearized static     &   1.02 &   3.31 \\ 
                                    & linearized dynamic    &   1.02 &   3.31 \\ 
                                    & velocity difference   &   1.02 &   2.75 \\ 
                                    & LAPACK {\tt DGEEV}    &   1.03 &   2.65 \\ 
            \hline 
            Bottom Layer Velocities & linearized static     &   1.65 &   2.72 \\ 
                                    & linearized dynamic    &   1.65 &   2.72 \\ 
                                    & velocity difference   &   1.33 &   2.82 \\ 
                                    & LAPACK {\tt DGEEV}    &   1.39 &   2.83 \\ 
            \hline \hline
        \end{tabular}
    \end{center}
    \caption{Calculated order of test cases run with a jump in bathymetry.
    \label{table:convergence_dry_wave}}
\end{table}

% subsubsection Eigensolver Comparisons and Convergence (end)

% subsection Simple Waves (end)

\subsection{Well-Balanced Test} % (fold)
\label{sub:well_balancing}

\begin{figure}[tbh]
    \centering
    % \begin{center}
        \includegraphics[width=0.42\textwidth]{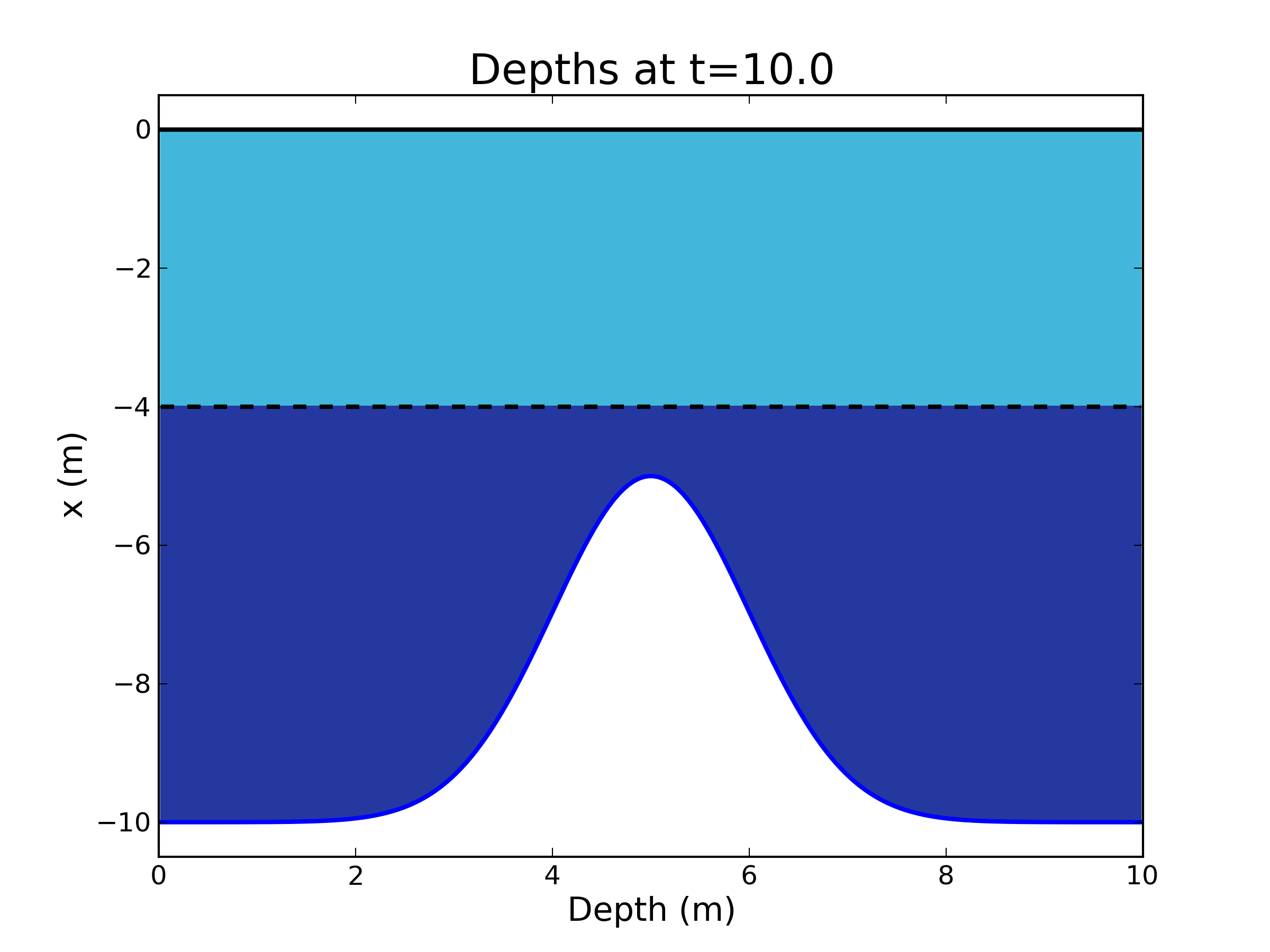}
        \includegraphics[width=0.42\textwidth]{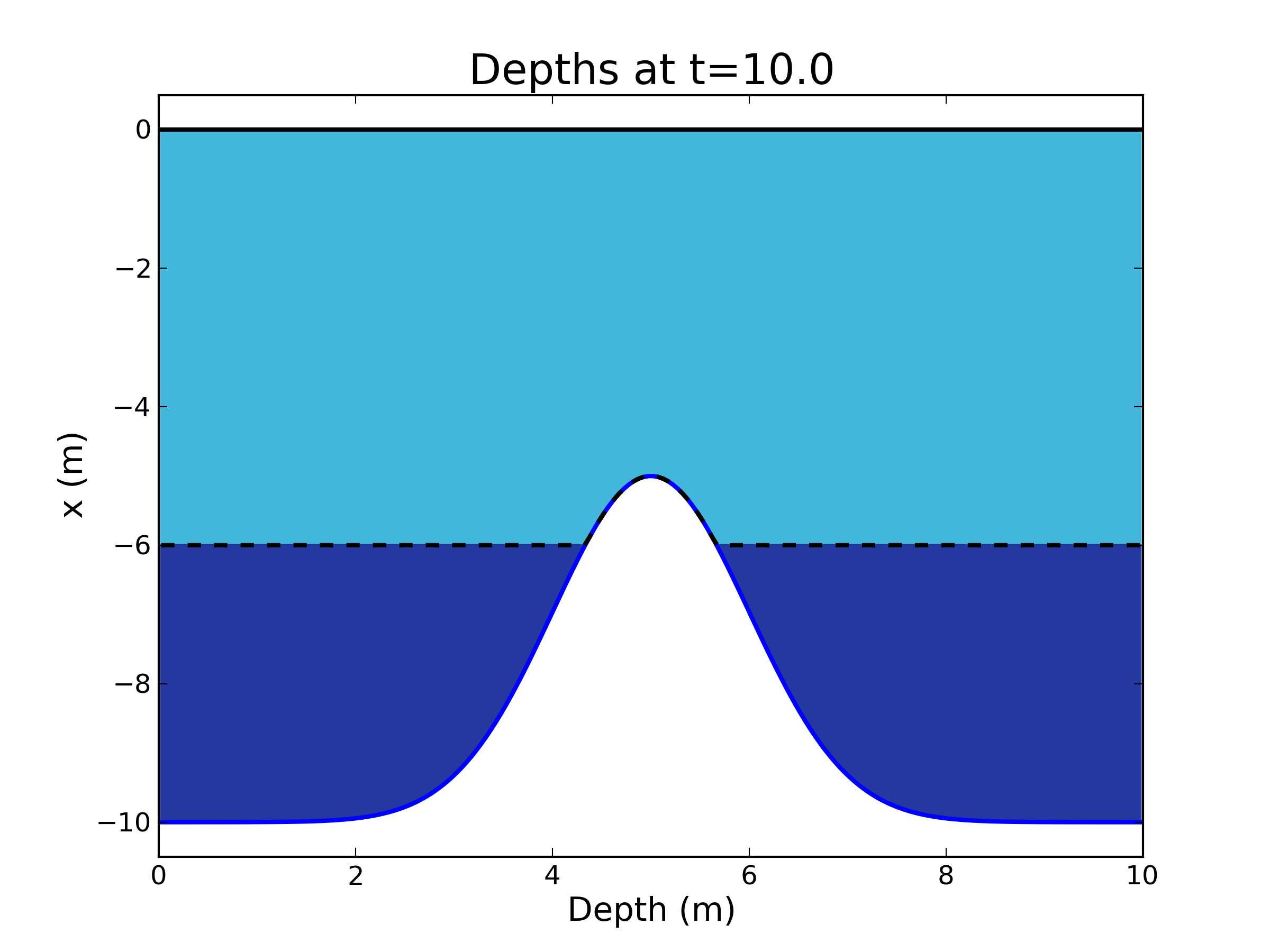} \\
        \includegraphics[width=0.42\textwidth]{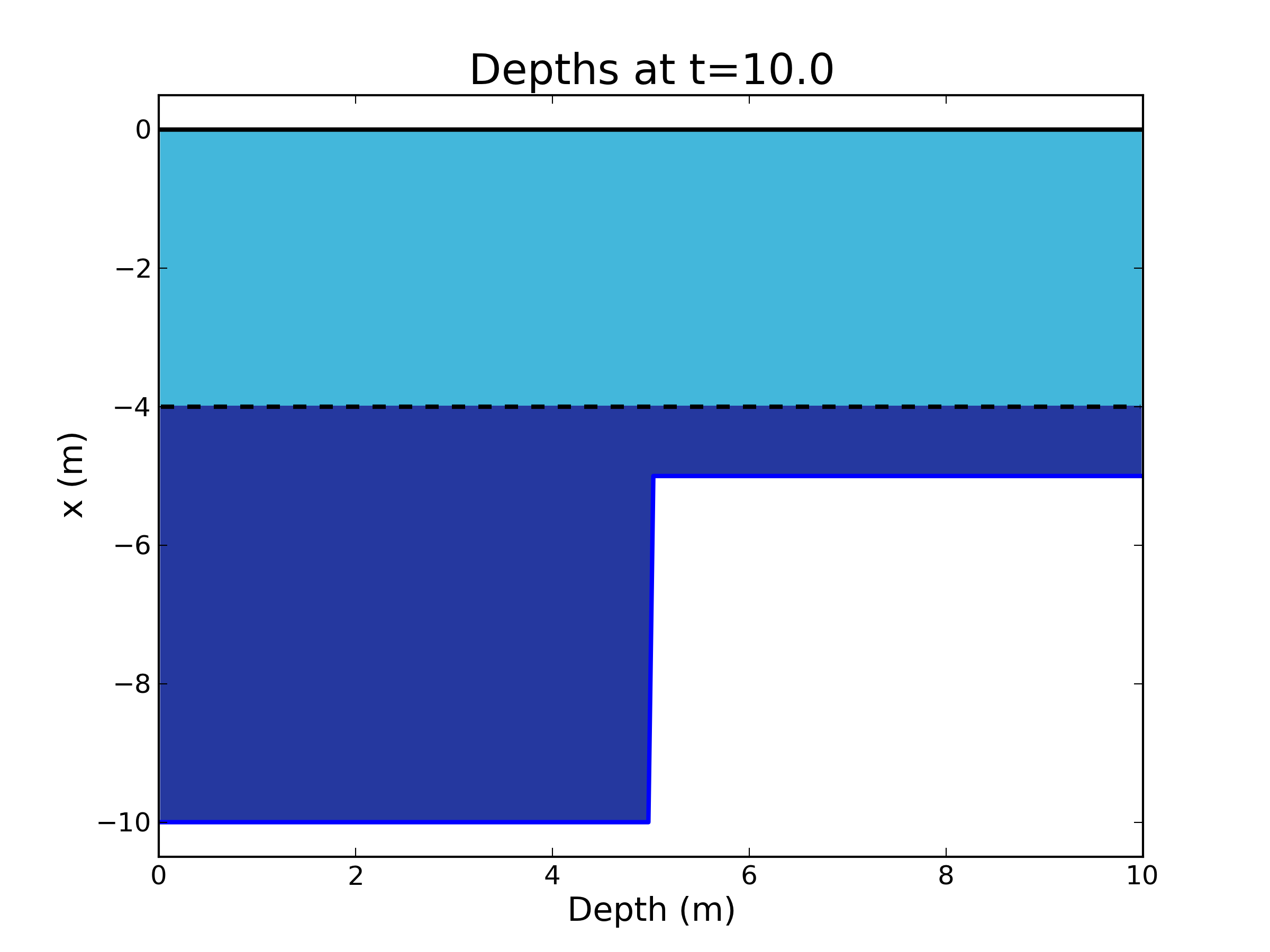}
        \includegraphics[width=0.42\textwidth]{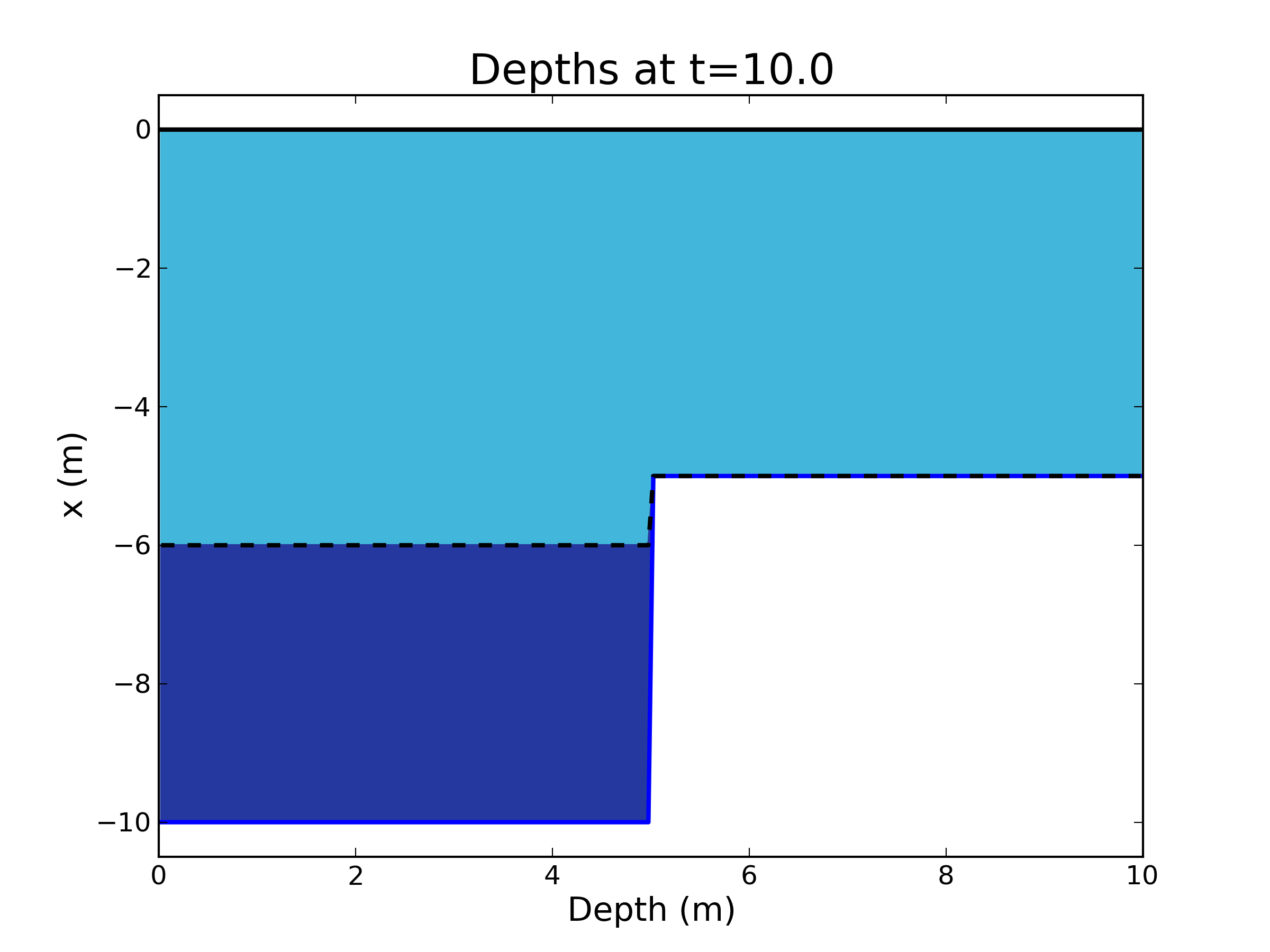}
    % \end{center}
    \caption{Illustration of the four well-balancing test cases after having been evolved to ten seconds.}
    \label{fig:well_balanced}
\end{figure}

Demonstration of the well-balanced properties of the solver in four cases was performed using two bathymetry configurations
\[
    b_\text{smooth} = -10 + 5 e^{-2 (x - 5)^2 / 5}
    ~~~~\text{and}~~~~ 
    b_\text{jump}(x) = \left \{ \begin{aligned}
        & -10  &\text{if}~x< 5& \\
        & -5 &\text{if}~x \ge 5
    \end{aligned} \right .
\]
illustrated in figure~\ref{fig:well_balanced} and two initial conditions with $\eta_1 = 0$ in both initial conditions and $\eta_2 = -4$ for the ``wet'' case and $\eta_2 = -6$ in the ``dry'' case.  The velocities of each layer are set to zero.  The final time for each simulation was $t=10$ seconds where the exact at rest solution was used to compute the $\text{L}^1$ errors in table~\ref{table:well_balanced_L_1} and $\text{L}^\infty$ in table~\ref{table:well_balanced_L_infty}.  Most error values agree with the precision used (double precision) but notably the depth and surfaces reduce to single precision in the smooth cases.  This reduction of accuracy may be due to conversion of the conserved quantities from $\rho_i h_i$ to $h_i$ and $\eta_i$ or possibly because of being near an inundation state.  

\begin{table}
    \begin{center}
        \begin{tabular}{||l|l|l|lll||}
        \hline \hline
        ~ & ~ & ~ & \multicolumn{3}{|c||}{$\text{L}^1$ Error} \\
        Bathymetry & Dry-State & Layer & $\rho_1 h_1$ & $\rho h u$ & $\eta$ \\
        \hline
        Smooth     &     False & 1 &               $0.00$ & $1.20\times10^{-11}$ & $8.88\times10^{-14}$ \\
                   &           & 2 & $8.26\times10^{-14}$ & $1.07\times10^{-11}$ & $7.19\times10^{-14}$ \\
        Smooth     &      True & 1 &  $1.17\times10^{-7}$ & $1.30\times10^{-11}$ & $1.62\times10^{-7}$  \\
                   &           & 2 &  $4.28\times10^{-8}$ & $5.28\times10^{-12}$ & $4.28\times10^{-8}$  \\
        Jump       &     False & 1 &               $0.00$ &               $0.00$ &              $0.00$  \\
                   &           & 2 &               $0.00$ &               $0.00$ &              $0.00$  \\ 
        Jump       &      True & 1 &               $0.00$ & $3.17\times10^{-11}$ & $8.88\times10^{-14}$ \\
                   &           & 2 & $1.48\times10^{-16}$ & $2.46\times10^{-12}$ &              $0.00$  \\
        \hline \hline
        \end{tabular}
    \end{center}
    \caption{$\text{L}^1$ errors computed from the true at-rest state.
    \label{table:well_balanced_L_1}}
\end{table}

\begin{table}
    \begin{center}
        \begin{tabular}{||l|l|l|lll||}
        \hline \hline
        ~ & ~ & ~ & \multicolumn{3}{|c||}{$\text{L}^\infty$ Error} \\
        Bathymetry & Dry-State & Layer & $\rho_1 h_1$ & $\rho h u$ & $\eta$ \\
        \hline
        Smooth     &     False & 1 &               $0.00$ & $9.20\times10^{-14}$ & $1.78\times10^{-15}$ \\
                   &           & 2 & $8.88\times10^{-16}$ & $7.95\times10^{-14}$ & $8.88\times10^{-16}$ \\
        Smooth     &      True & 1 &  $7.80\times10^{-9}$ & $1.24\times10^{-13}$ &  $7.96\times10^{-9}$ \\
                   &           & 2 &  $4.40\times10^{-9}$ & $6.01\times10^{-14}$ &  $4.40\times10^{-9}$ \\
        Jump       &     False & 1 &               $0.00$ &               $0.00$ &               $0.00$ \\
                   &           & 2 &               $0.00$ &               $0.00$ &               $0.00$ \\
        Jump       &      True & 1 &               $0.00$ & $2.27\times10^{-13}$ & $8.88\times10^{-16}$ \\
                   &           & 2 & $1.48\times10^{-16}$ & $6.78\times10^{-14}$ &               $0.00$ \\
        \hline \hline
        \end{tabular}
    \end{center}
    \caption{$\text{L}^\infty$ errors computed from the true at-rest state.
    \label{table:well_balanced_L_infty}}
\end{table}

% subsection well_balancing (end)

\subsection{Baroclinic Wave Wetting and Drying} % (fold)

\begin{figure}[tb]
    \centering
    % \begin{center}
        \includegraphics[width=0.42\textwidth]{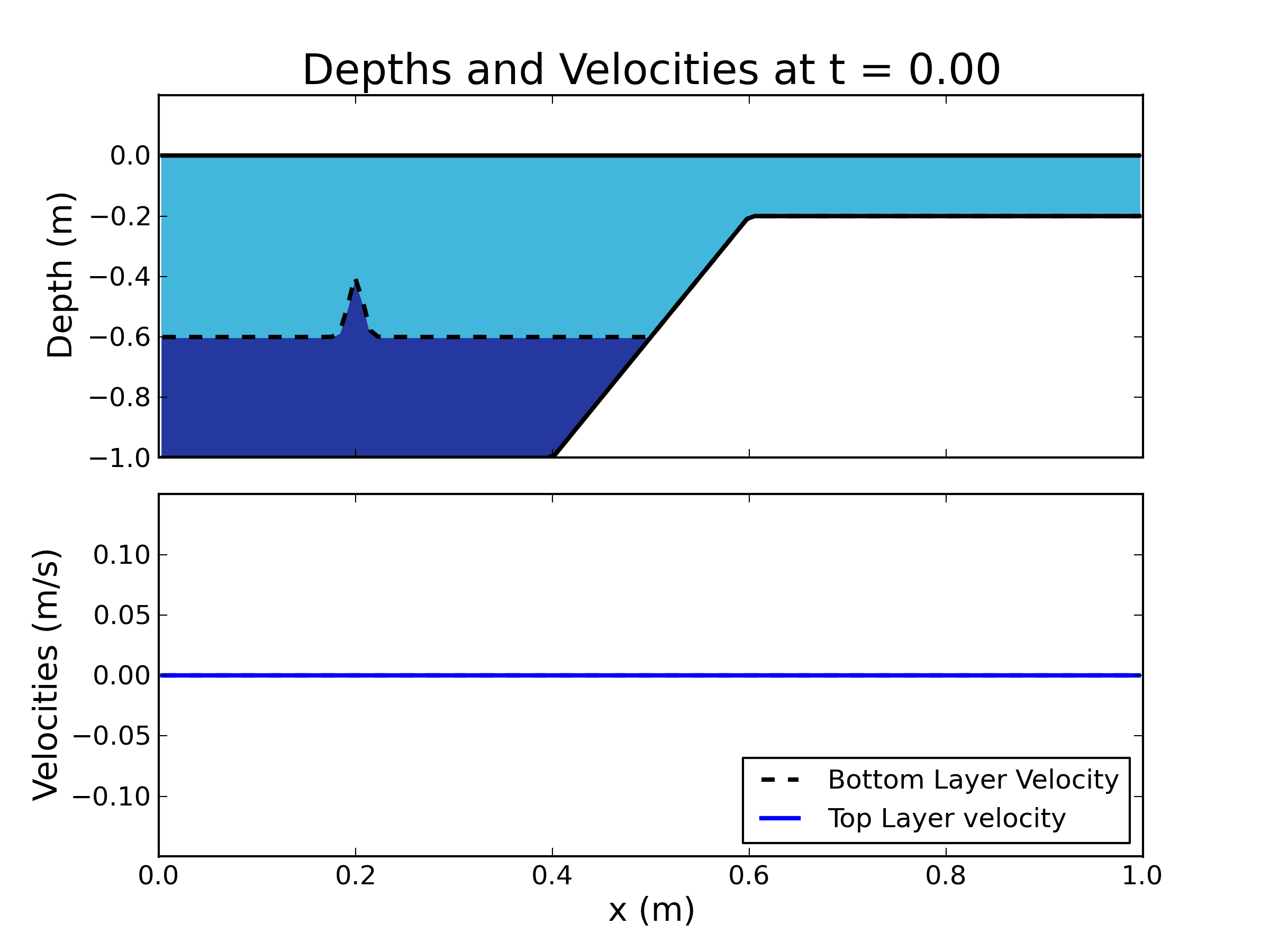}
        \includegraphics[width=0.42\textwidth]{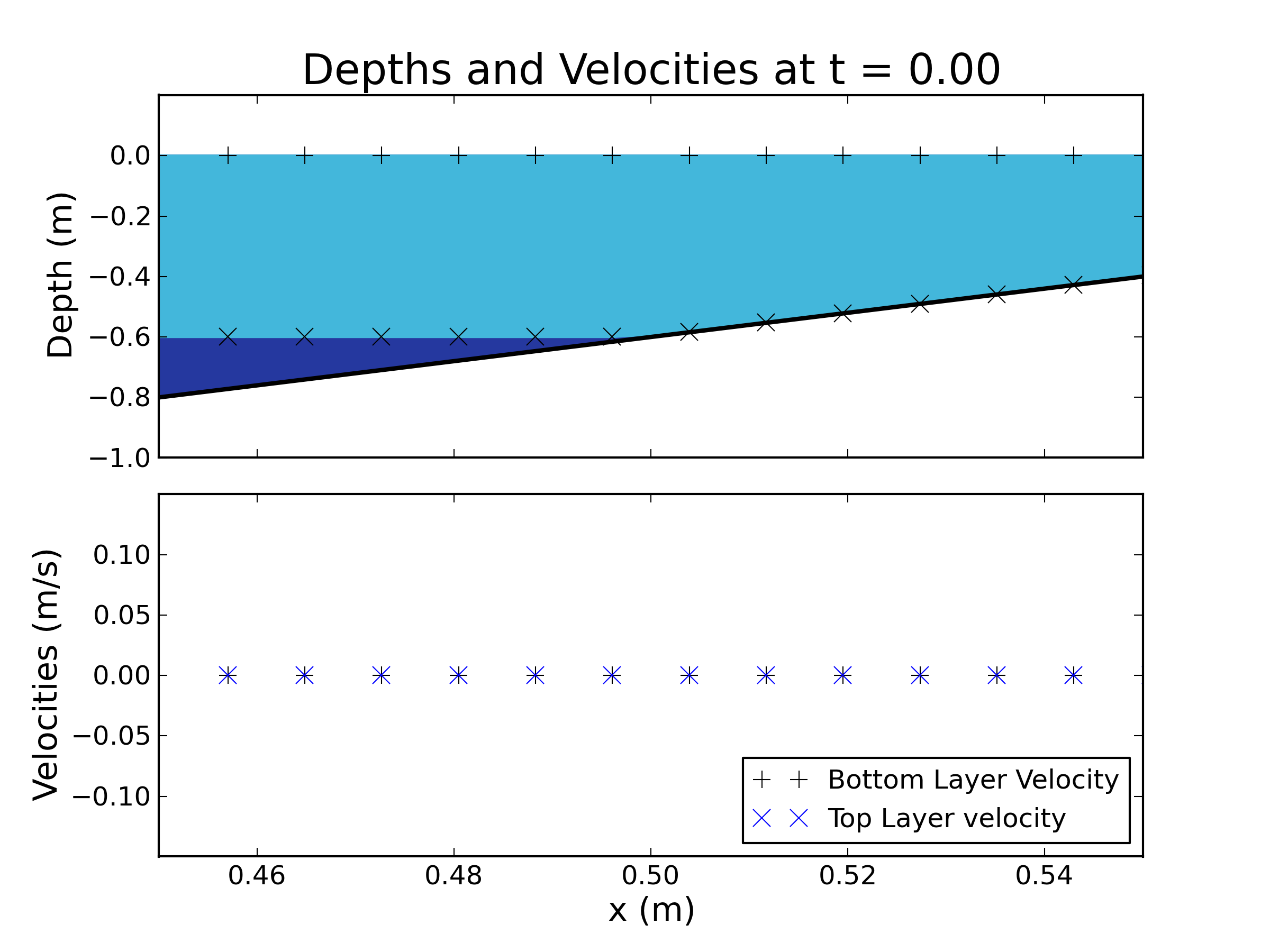} \\
        \includegraphics[width=0.42\textwidth]{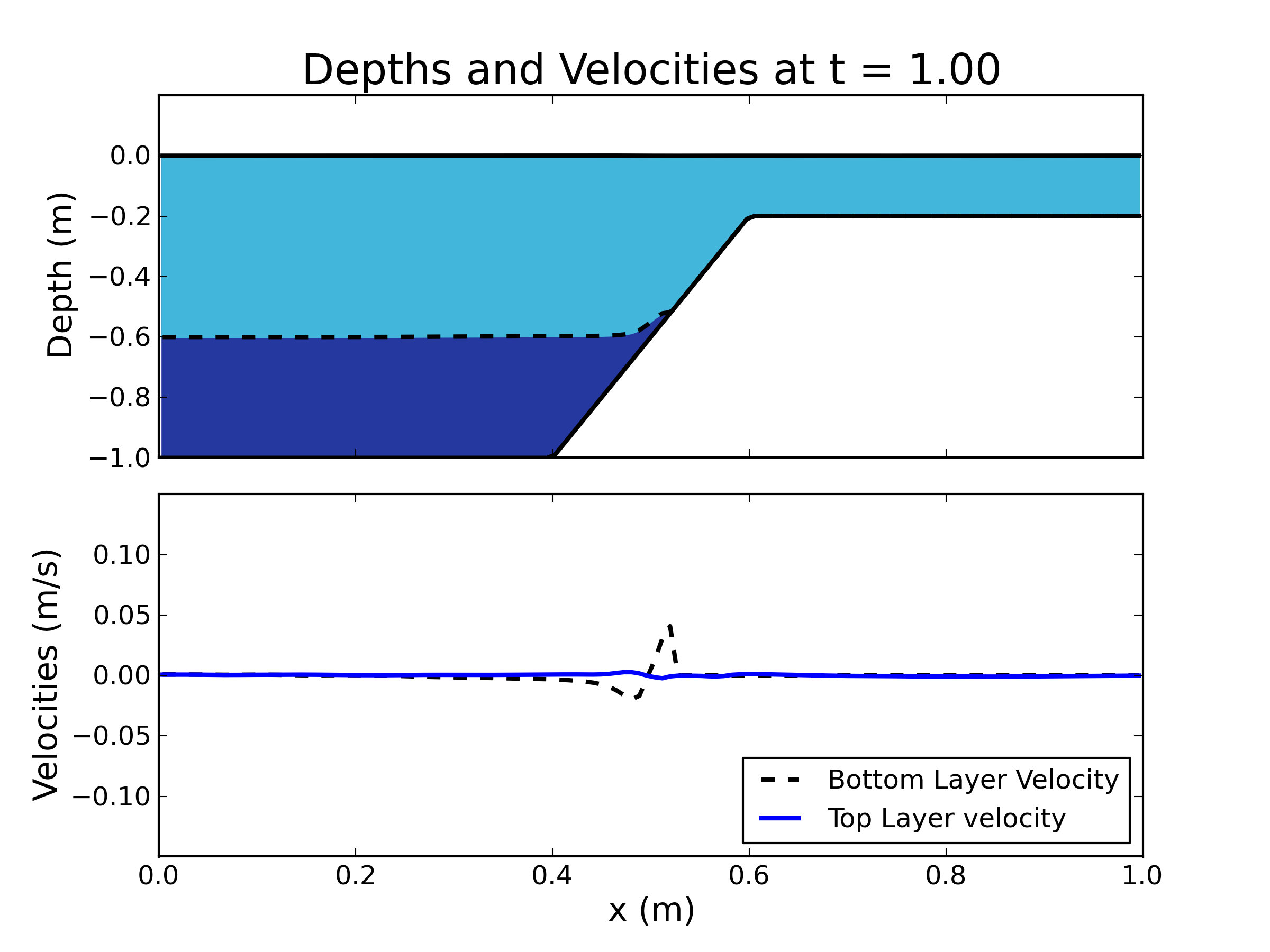}
        \includegraphics[width=0.42\textwidth]{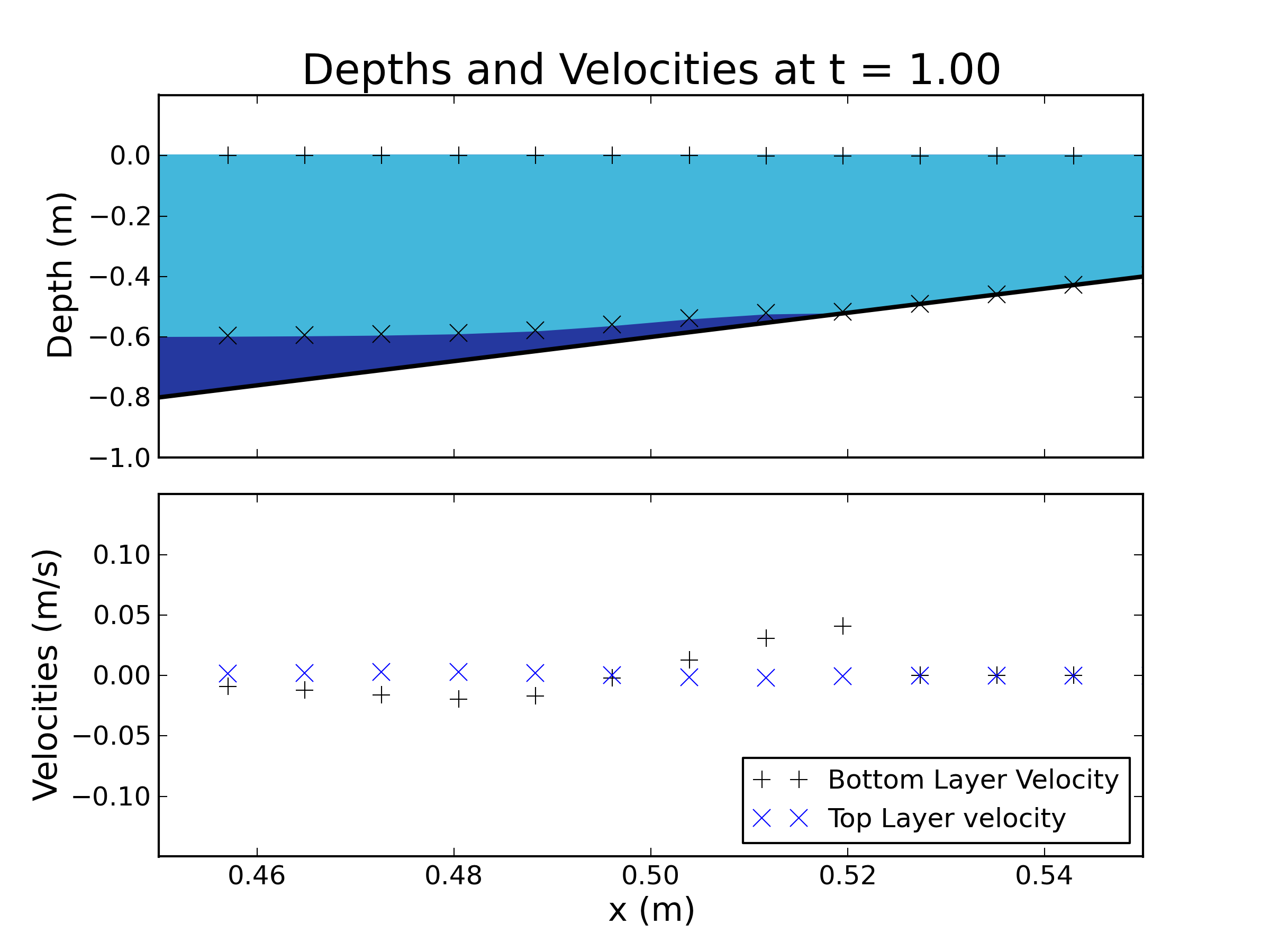}
    % \end{center}
    \caption{A Gaussian perturbation to the internal surface wets the slope located in the center of the domain.  These figures show the results of a simulation using the linearized dynamic eigensolver with 128 grid cells.  The left figure shows the entire domain and the right a zoomed-in section of the same domain located at the wet-dry interface of the bottom layer.}
    \label{fig:lapping_internal_wave}
\end{figure}

In order to demonstrate a case where wetting and drying occurs in the bottom layer, the example illustrated in figure~\ref{fig:lapping_internal_wave} was run using the methods outlined above for a variety of different resolutions.  The background state was set to $\hat{\eta}_1 = 0.0$ and $\hat{\eta}_2 = -0.6$ with $r=0.95$ and a Gaussian perturbation $\tilde{\eta}_2 = A \cdot e^{- ((x - x_0) / \sigma)^2}$ ($A=0.2$, $x_0 = 0.2$, and $\sigma=0.01$) added to the internal surface and allowed to propagate until the wave exited the domain.  The bathymetry was set to
\[
    b(x) = \left \{ \begin{aligned}
        &b_0  &\text{if}~x<x_0& \\
        & m \cdot (x - x_0) + b_0 &\text{if}~x \ge x_0 \le x_1 \\
        & b_1 &\text{if}~x \ge x_1
    \end{aligned} \right .
\]
where $b_0 = -1.0$, $x_0 = 0.4$, $b_1 = -0.2$, $x_1 = 0.6$, and $m = (b_0 - b_1) / (x_0 - x_1)$.  A Manning's-N type friction law was also applied to the bottom wet layer.  This term takes the form 
\[
    \rho_i h_i u_i^2 \frac{g n^2}{h_i^{4/3}}
\]
and is added via a first-order operator splitting approach with $n=0.022$.

Two interesting observations can be drawn from this example.  The first is that the method presented maintains positivity even in this case where small depths are present.  The second is that the ``true'' solution resolution used previously caused the simulation to become unstable.  This is probably the result of the simulation becoming non-hyperbolic near the dry-interface which is stabilized partially by numerical diffusive error present in lower resolutions.  Physically it is expected that overturning may occur at this interface and the diffusion may provide a means to handle the overturning numerically.  As mentioned previously, this loss of hyperbolicity is not seen in practice on the oceanic scales we are interested in representing.

% subsection Baroclinic Wetting-Drying (end)

\subsection{Idealized Ocean Shelf} % (fold)
\begin{figure}[ht] %  figure placement: here, top, bottom, or page
    \centering
    \subfloat[]{\includegraphics[width=0.38\textwidth]{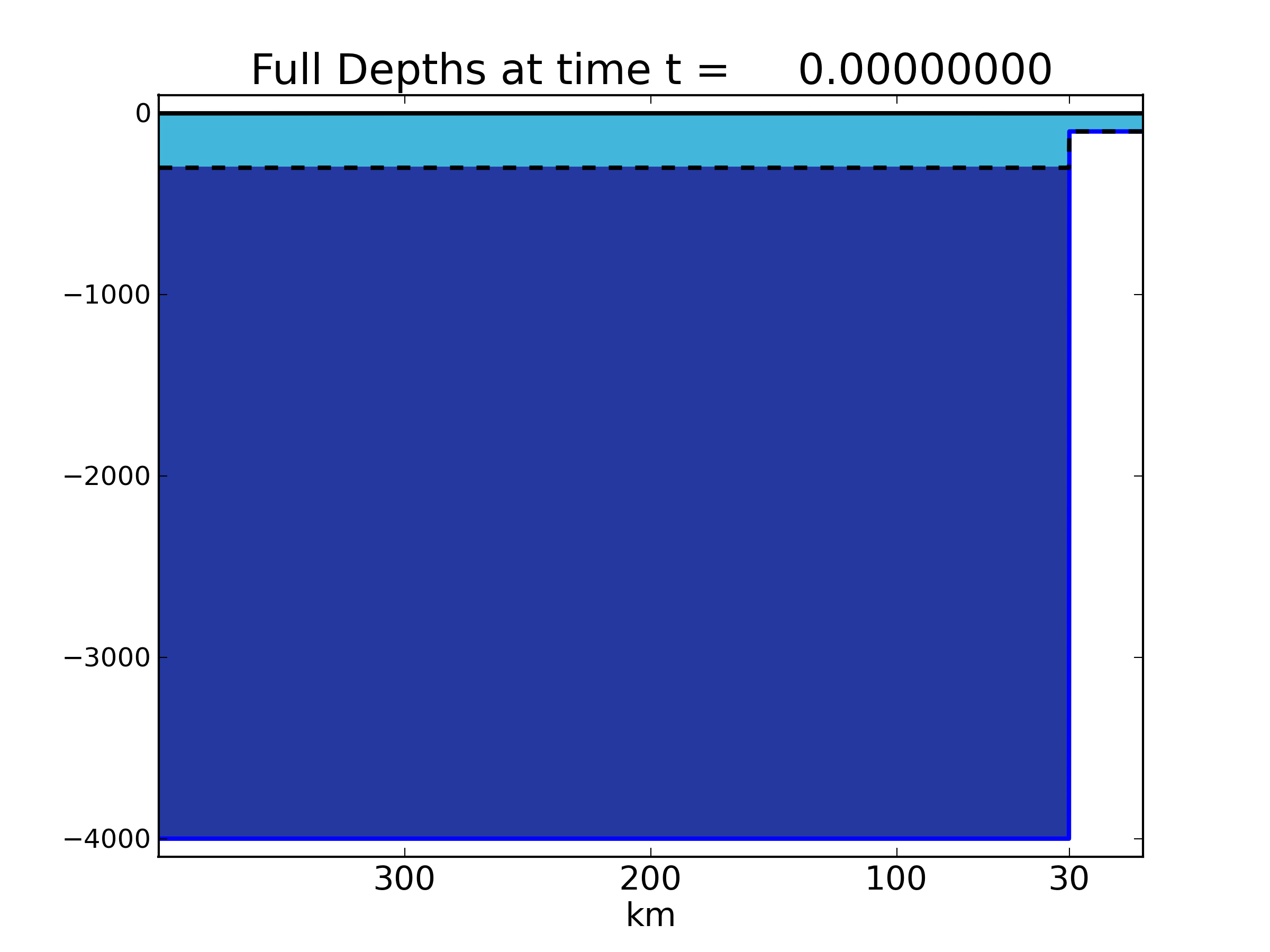} \label{fig:jump_shelf_bathy}}
    \subfloat[]{\includegraphics[width=0.30\textwidth]{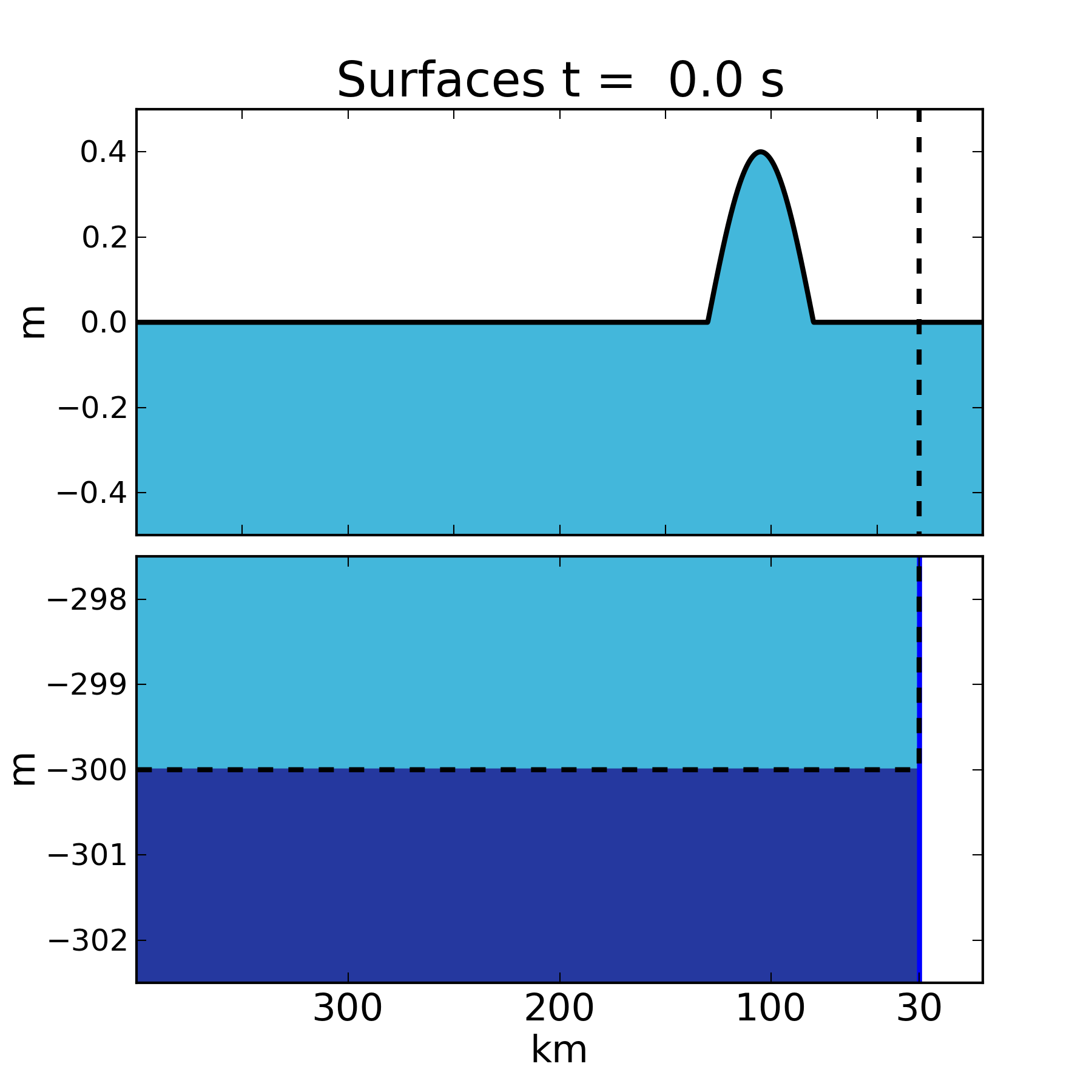} \label{fig:jump_shelf_init_surface}}
    % \subfloat[]{\includegraphics[width=0.35\textwidth]{jump_shelf/frame0000fig14.png} \label{fig:jump_shelf_init_velocties}}
    \caption{Problem setup for the idealized ocean shelf problem.  The vertical dotted line represents the location of the bathymetry jump.}
    \label{fig:jump_shelf_init}
\end{figure}

\begin{figure}[htpb] %  figure placement: here, top, bottom, or page
    \centering
    \includegraphics[width=0.3\textwidth]{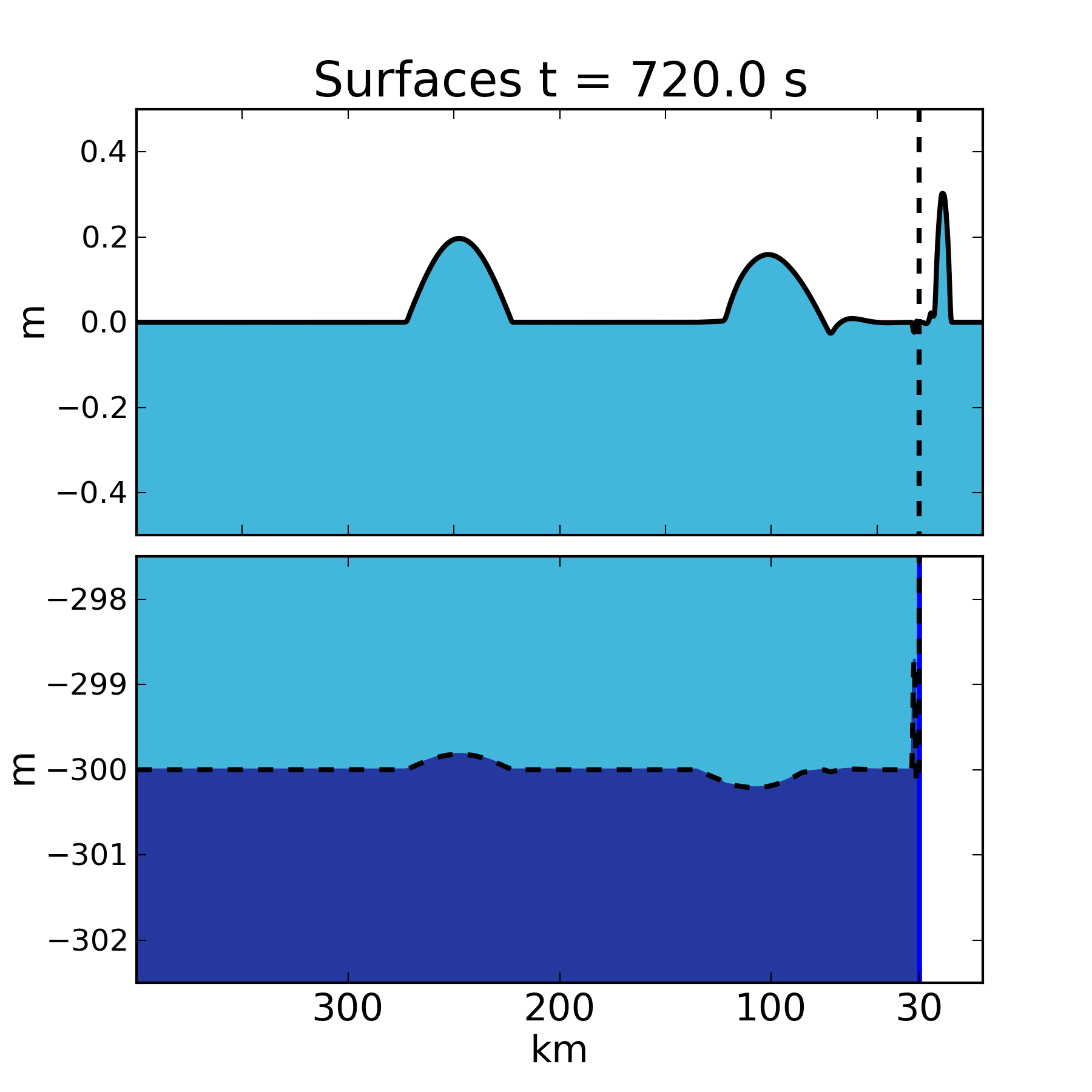}   
    \includegraphics[width=0.3\textwidth]{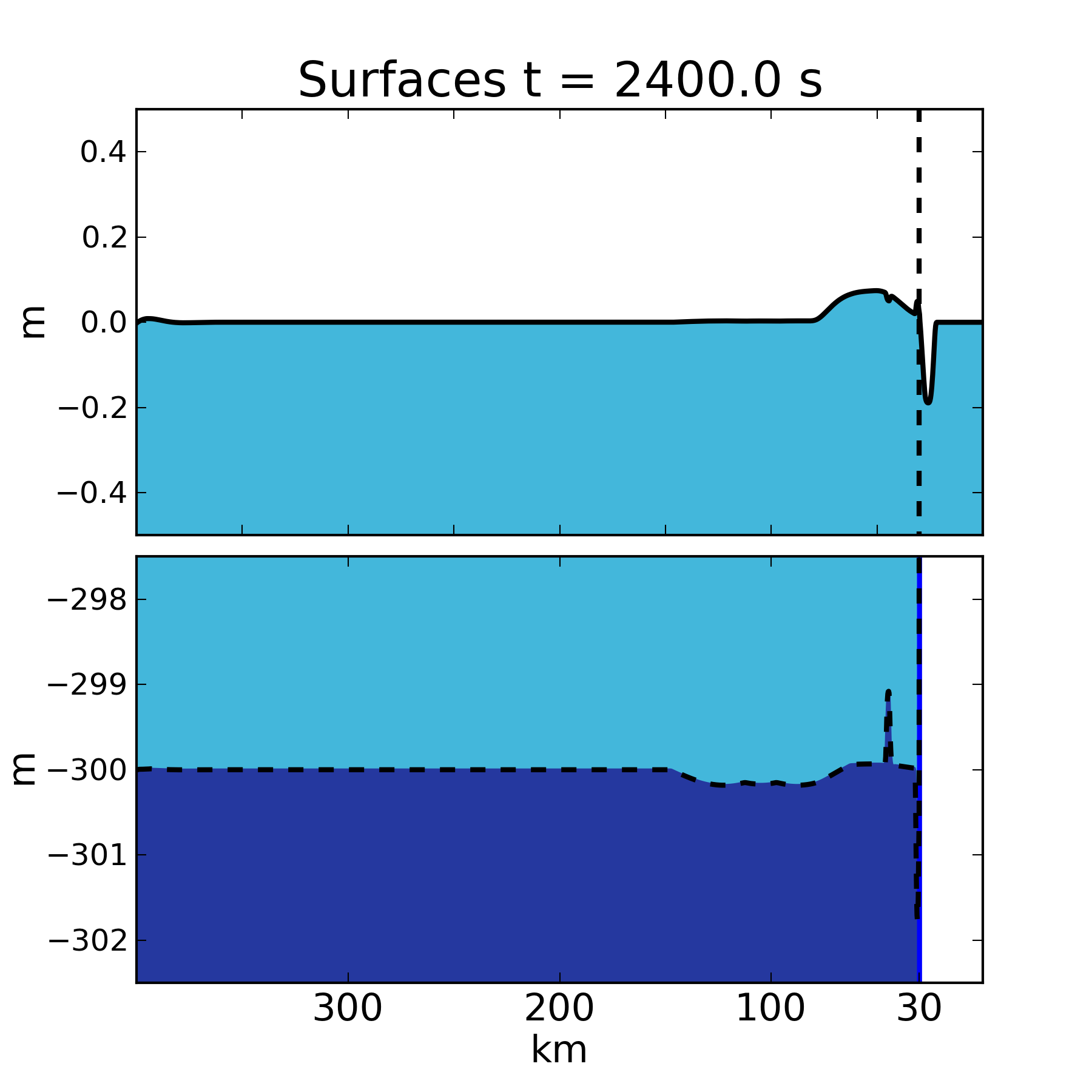}   
    \includegraphics[width=0.3\textwidth]{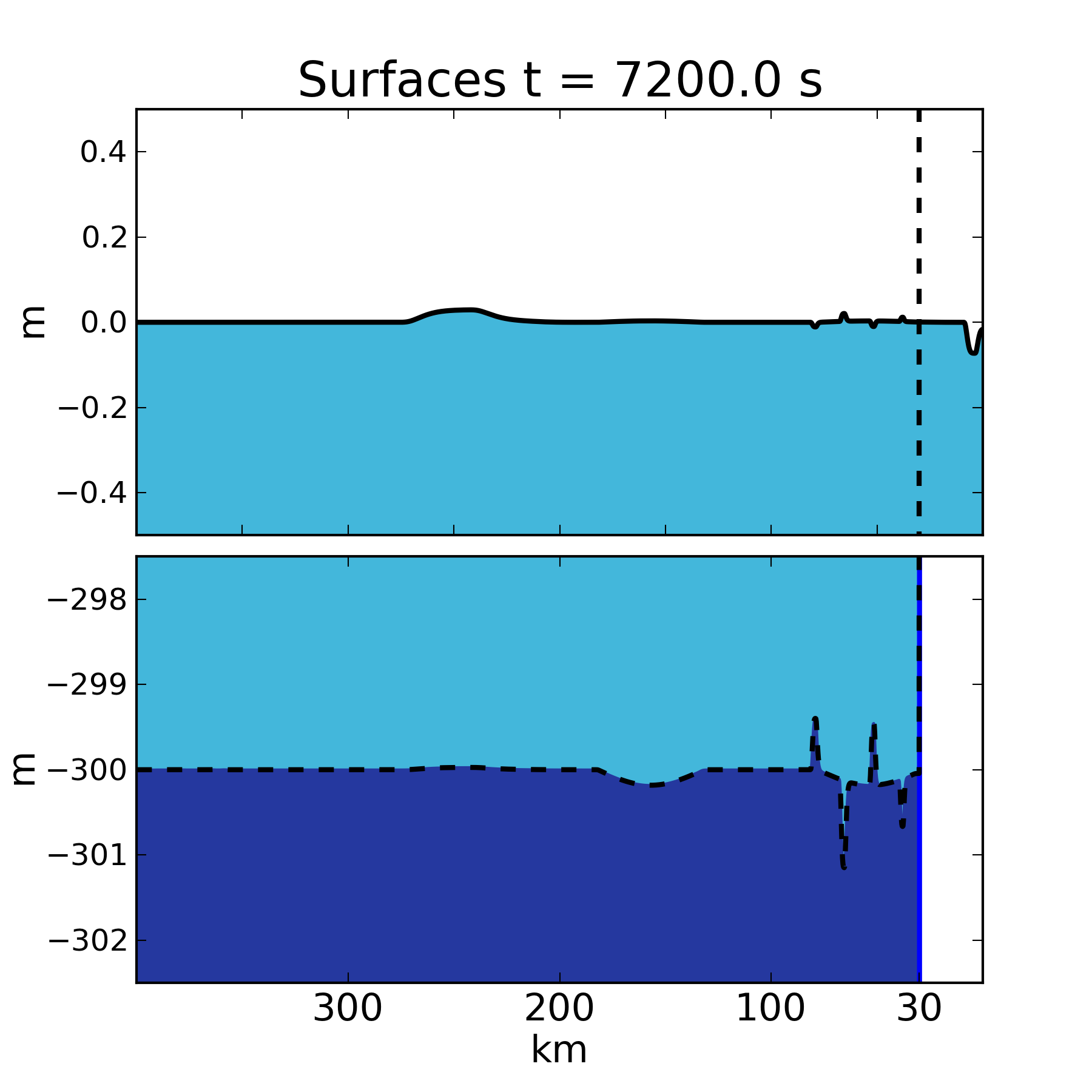} \\
    \includegraphics[width=0.3\textwidth]{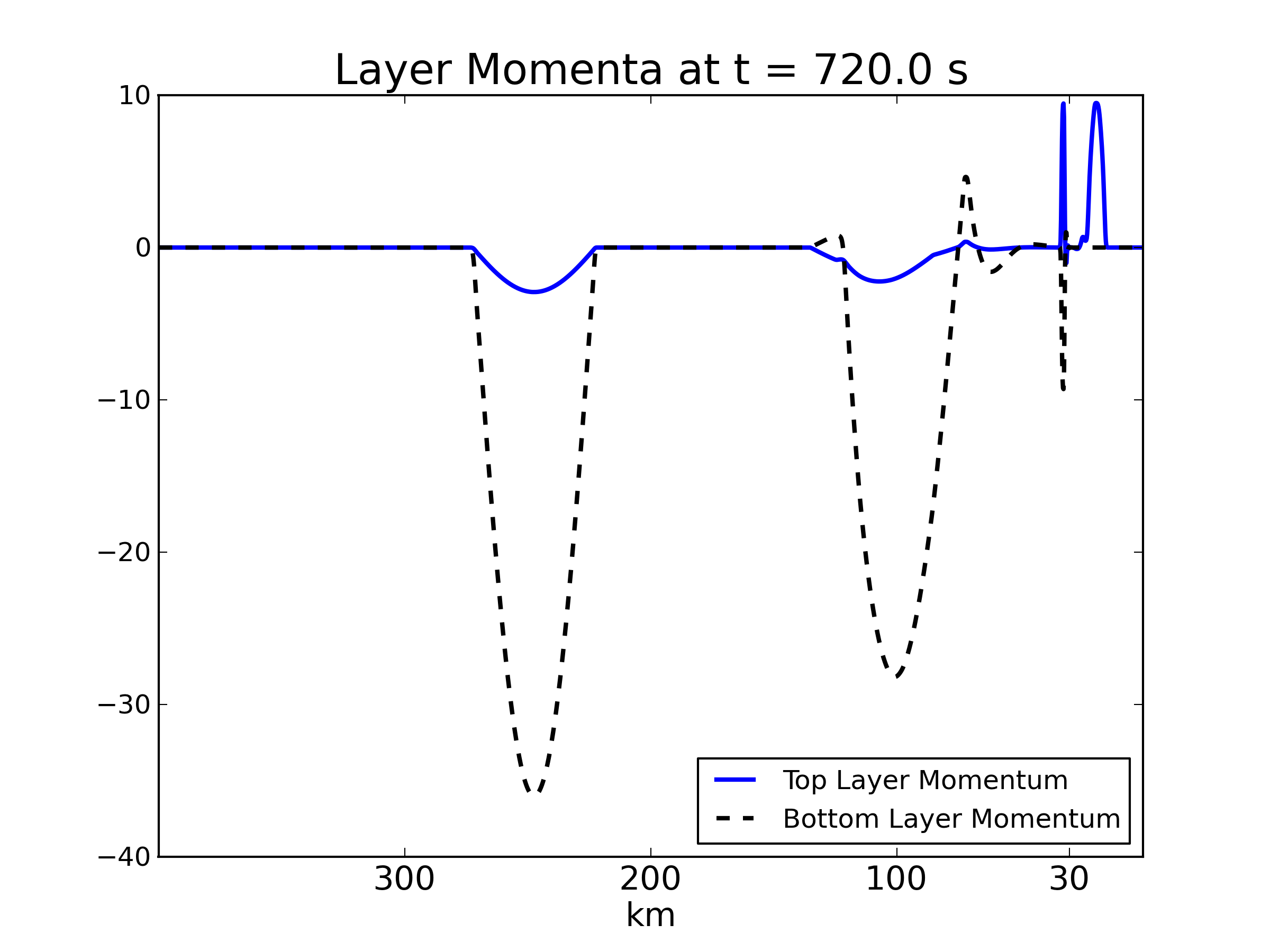}
    \includegraphics[width=0.3\textwidth]{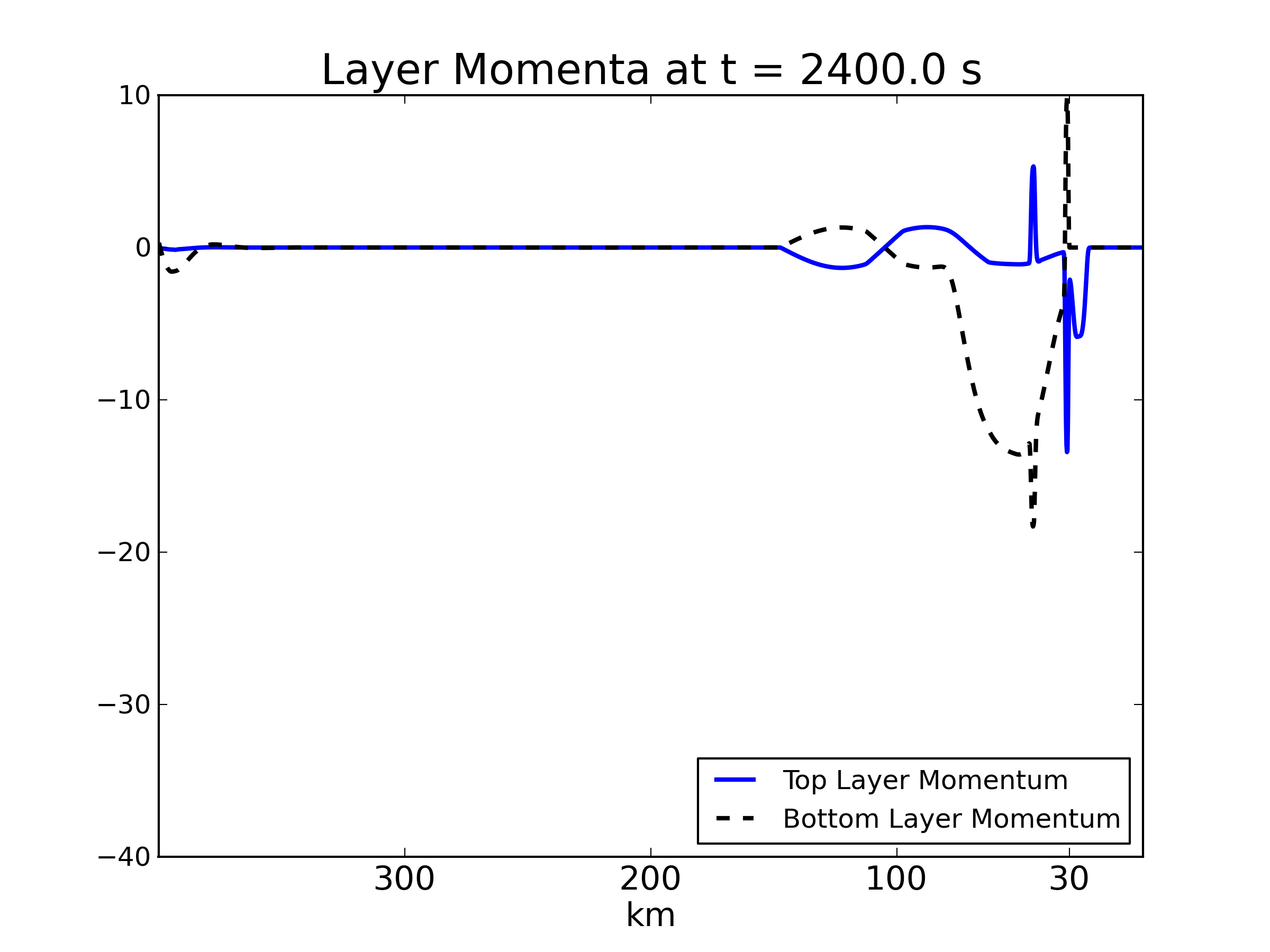}
    \includegraphics[width=0.3\textwidth]{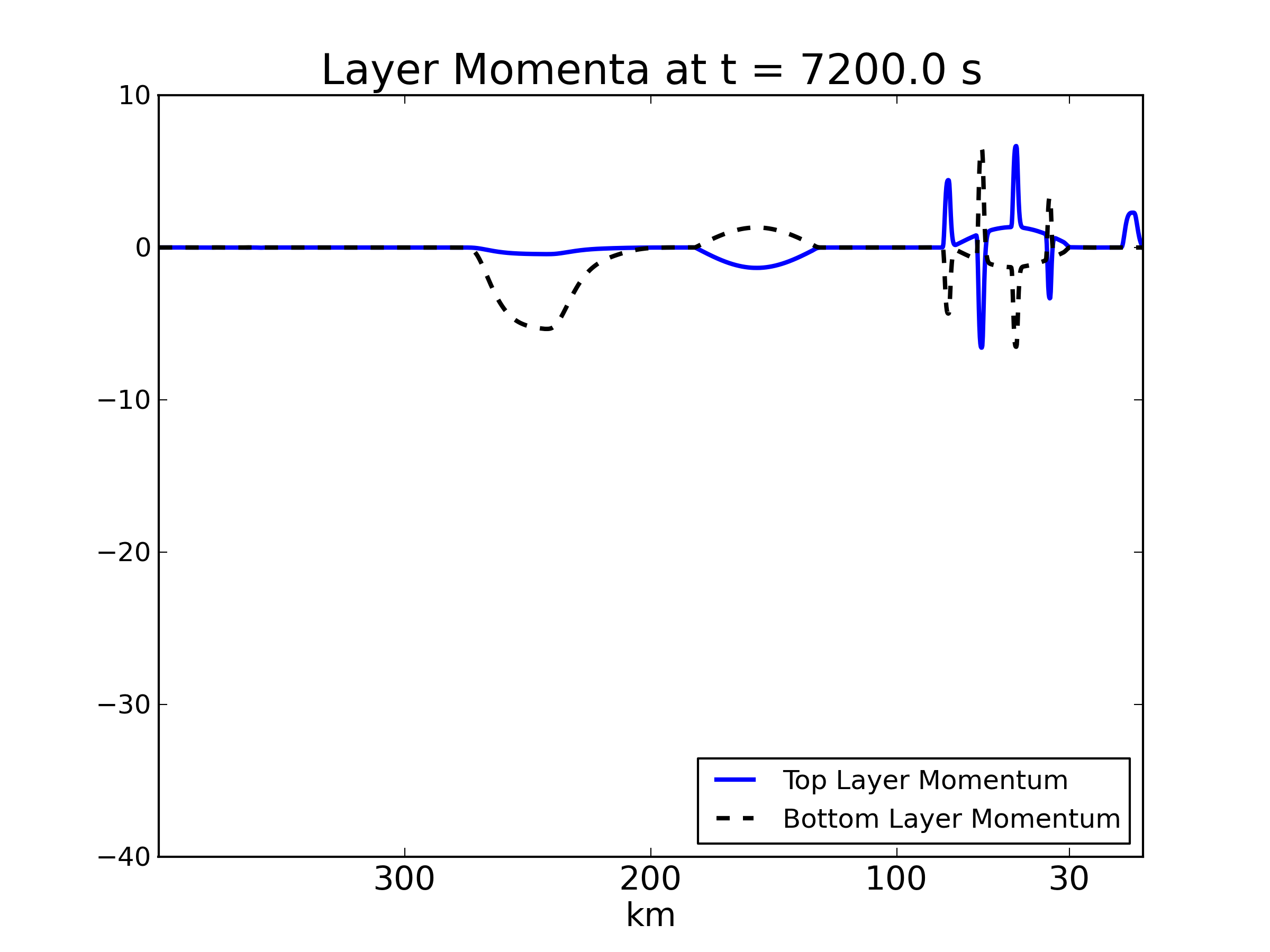}
    \caption{Time snapshots from an idealized ocean shelf test using the dynamic linearized eigensolver.}
    \label{fig:jump_shelf_solution}
\end{figure}

\begin{figure}[htb] %  figure placement: here, top, bottom, or page
    \centering
    \includegraphics[width=\textwidth]{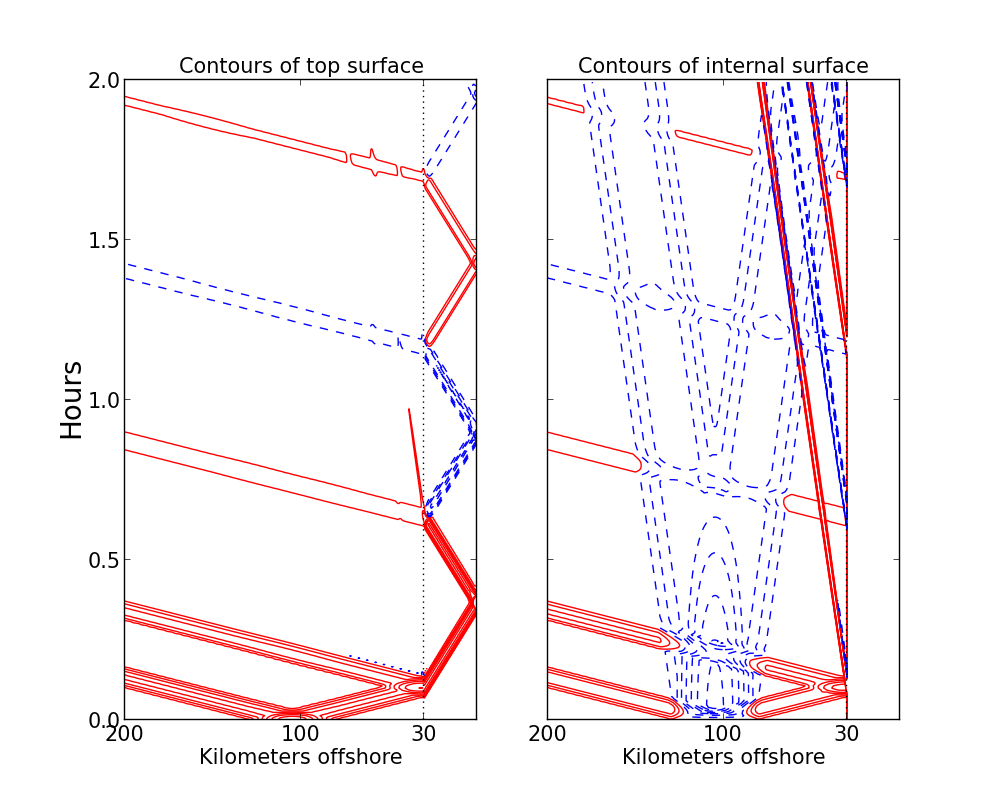} 
    \caption{Contour plots of the top and internal surface height through time for the oceanic test.  The top surface has 30 contour levels from -0.4 to 0.4 and bottom surface has 30 contour levels from -0.5 to 0.5 where solid contours indicate positive displacement and dashed contours negative displacement.  The dotted line is a reference for the bathymetry changes.}
    \label{fig:jump_shelf_contour}
\end{figure}

\label{sub:ocean_shelf}
As a test of the solver in a context similar to a large oceanic basin, we will look at a test originally proposed in \cite{Berger:2011du} for the single-layer equations.  It involves a large domain size with a step discontinuity representing an idealized continental shelf and a wall boundary condition at the left boundary representing a coast.  The initial condition is a perturbation approximately $0.4~\text{m}$ in amplitude of both the sea and internal surface of an ocean at rest.  The densities of the layers are $\rho_1 = 1025~\text{kg} / \text{m}^3$ and $\rho_2 = 1045~\text{kg} / \text{m}^3$, similar to that found in the ocean.  Figure~\ref{fig:jump_shelf_bathy} details the bathymetry where $b_\ell = -4000$ and $b_r = -100$ with the shelf 30 kilometers from the right boundary and an ocean at rest specified by $\hat{\eta}_1 = 0$ and $\hat{\eta}_2 = -300$.  The perturbation to the surface is 
\[
    \tilde{\eta}_1(x) = \epsilon \sin \left[\frac{\pi (x - x_\text{mid})}{-80~\text{km} - x_\text{mid}}\right]
\] 
where $\epsilon = 0.4$ and $x_\text{mid} = -130~\text{km}$.  For the simulations depicted 2000 grid cells were used.  No friction terms were added in this example.  Figure~\ref{fig:jump_shelf_solution} shows snapshots of the computed solution.  Note that the x-axis has been flipped in the plots to match what was done in \cite{Berger:2011du}.  Also following \cite{Berger:2011du} contour plots of the top and internal surface can be found in figure~\ref{fig:jump_shelf_contour}.

An interesting feature of this solution are the high amplitude and short wavelength internal waves being generated as the waves are coming off of the shelf.  This is due to the bottom layer reacting to the immediate loading of mass from the shelf.  Since the internal wave speeds are significantly smaller than the shallow water gravity speeds from the shelf, the waves have much shorter wavelengths as can be seen in figure~\ref{fig:jump_shelf_contour} where the internal surface contours show two significant wave speeds.  Simulations were also performed at multiple grid resolutions and for multiple eigenspace solvers which all indicated that these high amplitude and short wavelength internal waves are present in the two layer model and not a numerical instability.

% subsection idealized_ocean_shelf (end)

% section results (end)
% ============================================================================

% ============================================================================
\section{Concluding Remarks} \label{sec:conclusions} % (fold)

We have developed a numerical method for solving the multilayer shallow water equations utilizing an f-wave propagation finite volume method.  The approach robustly handles dry-states and allows for the use of existing single layer solvers when applicable.  A number of eigensolvers have been shown to perform reasonably well except at dry-states where the linearized eigensolver can be used.  The solver also maintains the nonlinear behavior expected in the system of equations even though a linearization of the eigenspace has been used.

This approach appears promising for a variety of applications to oceanic flows that may have simple vertical structure or forcing that is vertically localized such as wind or friction.  Part of this promise is the possible leveraging of current state-of-the-art methods available to single-layer shallow water methods, such as adaptive mesh refinement, on this class of applications.  To this end, the method has been implemented in the GeoClaw software package  \cite{Berger:2011vi} with future work needed in the algorithms for refining and coarsening.  Application studies are also under way, in particular for storm surge simulations, exploring the benefit of the increased vertical structure over single-layer approaches.  Further extensions of the method to correct for the entropy violations and possible refinement of the inundation wave speeds are also of particular interest.

% section conclusions (end)
% ============================================================================

\vskip 10pt
{\bf Acknowledgments.}
The author would like to thank Randall LeVeque, Kristen Thyng and the reviewers for their time and effort in editing this article.  This research was supported in part by DOE grant DE-FG02-88ER25053, NSF grants RTG-0838212 and DMS-0914942, and ONR Grant N00014-09-1-0649.

\bibliographystyle{elsarticle-num}
\bibliography{database}

\begin{thebibliography}{10}
\expandafter\ifx\csname url\endcsname\relax
  \def\url#1{\texttt{#1}}\fi
\expandafter\ifx\csname urlprefix\endcsname\relax\def\urlprefix{URL }\fi
\expandafter\ifx\csname href\endcsname\relax
  \def\href#1#2{#2} \def\path#1{#1}\fi

\bibitem{Salmon:2002vg}
R.~Salmon, {Numerical solution of the two-layer shallow water equations with
  bottom topography}, Journal of Marine Research 60 (2002) 605--638.

\bibitem{abgrall:2007}
R.~Abgrall, S.~Karni, {Two-Layer Shallow Water Systems: A Relaxation Approach},
  SIAM Journal of Scientific Computing 31~(3) (2009) 1603--1627.

\bibitem{Bouchut:2008kn}
F.~Bouchut, T.~Morales~de Luna, {An entropy satisfying scheme for two-layer
  shallow water equations with uncoupled treatment}, ESAIM: M2AN 42~(4) (2008)
  683--698.

\bibitem{CastroDiaz:2010kn}
M.~J. Castro-D{\'\i}az, E.~D. Fern{\'a}ndez-Nieto, J.~M. Gonz{\'a}lez-Vida,
  C.~Par{\'e}s-Madr{\~n}al, {Numerical Treatment of the Loss of Hyperbolicity
  of the Two-Layer Shallow-Water System}, Journal of Scientific Computing
  (2010) 1--29.

\bibitem{Castro:2001td}
M.~Castro, J.~Mac{\'\i}as, C.~Par{\'e}s, {A Q-scheme for a class of systems of
  coupled conservation laws with source term. Application to a two-layer 1-D
  shallow water system}, Esaim-Mathematical Modelling and Numerical
  Analysis-Modelisation Mathematique Et Analyse Numerique 35~(1) (2001)
  107--127.

\bibitem{LAWRENCE:1990ul}
G.~A. Lawrence, {On the hydraulics of Boussinesq and non-Boussinesq two-layer
  flows}, Journal of Fluid Mechanics 215 (1990) 457--480.

\bibitem{Bale:2002}
D.~S. Bale, R.~LeVeque, S.~Mitran, J.~A. Rossmanith, {A Wave Propagation Method
  for Conservation Laws and Balance Laws with Spatially Varying Flux
  Functions}, SIAM Journal of Scientific Computing 24~(3) (2002) 955--978.

\bibitem{LeVeque:2001}
R.~LeVeque, M.~Pelanti, {A Class of Approximate Riemann Solvers and Their
  Relation to Relaxation Schemes}, Journal of Computational Physics 172 (2001)
  572--591.

\bibitem{LeVeque:2002aa}
R.~LeVeque, {Finite Volume Methods for Hyperbolic Problems}, Cambridge Texts in
  Applied Mathematics, Cambridge University Press, Cambridge, UK, 2002.

\bibitem{Schijf:1953vz}
J.~B. Schijf, J.~C. Sch{\"o}nfled, {Theoretical considerations on the motion of
  salt and fresh water}, in: Proc. of the Minn. Int. Hydraulics Conv., Joint
  meeting IAHR and Hyd. Div ASCE., IAHR, 1953, pp. 321--333.

\bibitem{Mandli:2011te}
K.~T. Mandli, {Finite Volume Methods for the Multilayer Shallow Water Equations
  with Applications to Storm Surges}, Ph.D. thesis, University of Washington
  (Jun. 2011).

\bibitem{George:2008aa}
D.~L. George, {Augmented Riemann solvers for the shallow water equations over
  variable topography with steady states and inundation}, Journal of
  Computational Physics 227~(6) (2008) 3089--3113.

\bibitem{Berger:2011du}
R.~J. LeVeque, D.~L. George, M.~J. Berger, {Tsunami Propagation and inundation
  with adaptively refined finite volume methods}, Acta Numerica (2011)
  211--289.

\bibitem{Berger:2011vi}
M.~J. Berger, D.~L. George, R.~J. LeVeque, K.~T. Mandli, {The GeoClaw software
  for depth-averaged flows with adaptive refinement}, Advances in Water
  Resources 34~(9) (2011) 1195--1206.

\end{thebibliography}
\end{document}